# A Short-term Planning Framework for the Operation of Tanker based Water Distribution System in Urban Areas

**Abhilasha Maheshwari**[†#], **Shamik Misra**[†#]**, Ravindra D. Gudi**[†,*]**,**
**Senthilmurugan Subbiah**[‡]

[†] Department of Chemical Engineering, Indian Institute of Technology Bombay, Powai, Mumbai - 400076, India

[‡] Department of Chemical Engineering, Indian Institute of Technology Guwahati, Guwahati, Assam-781039, India

[#]Both authors have contributed equally.

*To whom correspondence must be addressed: Ravindra D. Gudi (e-mail: ravigudi@iitb.ac.in ,

Department of Chemical Engineering, IIT Bombay, Mumbai, India, contact: +91-022-2576-7231)

# Abstract

Tanker based distribution systems have been prevalent in developing countries to supply clean and pure water in different regions. To efficiently operate such tanker service system, a large fleet of tanker trucks are required to transport water among several water sources, water treatment plants and consumers spanning across the regions. This requires tighter coordination between water suppliers, treatment plant operations, and user groups to use available water resources in a sustainable manner, along with the assurance of water quality and timely delivery. This paper proposes a novel formulation to assist decision making for optimizing tanker based water distribution system and treatment operations, with an overall objective of minimizing total operating cost, such that all the constraints related to the water demand, supply operations, environmental and social aspects are honored, while supplying water to a maximum number of users. The problem is formulated and solved as an MILP and captures all the nuances related to (i) water availability limitations and quality constraints from different sources, (ii) maintaining water quality as it transports via tankers, (iii) water demands for various end-use purposes, and (iv) transportation across water supply chain. The proposed novel framework is applied to a realistic urban model to find the optimal tanker delivery schedule, ensuring appropriate treatment and timely delivery of water. The results of the case study conducted on a representative scale problem also elucidate several aspects of treatment plant operation and consumer demand fulfillment for efficient planning and management of tanker based water distribution system.

## 1. Introduction

Providing safe and adequate water supply to the growing population in rapidly expanding urban areas is one of the goals set by the United Nations for sustainable development.[1] The urban water distribution system is broadly operated in two ways: (i) piped water network and (ii) mobile vendors such as tanker truck vehicles.[2] Worldwide, it is considered a priority action to extend piped distribution network towards achieving 100 % connectivity and 24x7 supply of treated water. However, full coverage of the expanding urban regions with the piped network is expected to take a long time and involves challenges related to cross-contamination of water due to leakages, corrosion of pipes, etc.[3]

Even among the population which is currently served by the existing piped network in Indian metropolitan cities like Delhi, Mumbai, Hyderabad, and Chennai, less than 55% of the people receive water for only 3 hours in a day.[3] Furthermore, insufficient volumes at low pressure in such intermittent water supply methods fail to fully satisfy the consumer demands in urban areas.[4] The water quality from the intermittently operated pipe networks also gets affected due to the intrusion of contaminant caused by (1) negative pressure differentials in pipelines, (2) dissolution of corrosion scales from old pipes, and (3) microbial detachment from the pipe surfaces, etc.[5] Also, changing climate conditions resulting in frequent droughts and depleting groundwater levels have further increased the water demand-supply gap in urban cities.[6, 7]

The 2016 Safe Water Network Report for Mumbai city shows that Bombay Municipal Corporation managed to supply only 3200 MLD water through piped network against the city demand of 3856 MLD in the year 2011 and according to the water demand projections this demand-supply gap will increase to 1100 MLD by the year 2021.[8] Thus, the limitations in (i) water distribution from intermittently operated piped network, (ii) direct water use from arsenic contaminated groundwater sources[9], and (iii) huge seasonal variability in available water resources, have increased the reliance of people on tanker based water supply service as a quick and reliable alternative to fill the supply gap.[10] Statistics from a case study in Chennai, a major metropolitan city in India, show that tanker based water supply constitutes about 21 % of the total water supplied by the Chennai Metro Water Supply Board.[11] A research study by Londhe et al.[12] spanned over six cities having population over 1 lakh (class I cities of India), found that the large water consumers such as commercial establishments, hospitals, schools, etc. rely on tanker based water suppliers for their entire or at least partial water demands.

Several other research studies on households in cities like Mumbai, Hyderabad, and Delhi found the significant prevalence of tanker based water supply to fulfill public water needs throughout the year.[13, 14] Motivated by such factors, the Delhi Jal Board recently established a water supply system to serve potable water through 1000 tankers in eight zones across the Delhi City.[15] Such tanker based water distribution is also found to be prominent in recent water consumption simulation and survey case studies on cities in Jordan, Nigeria, Kenya, and other parts of the World.[16-18]

While the tanker based water distribution is considered an effective alternative to extend water supply service, aspects such as timely delivery and quality of water delivered to the consumer are crucial.[19] Srinivasan et al.[20] found in a survey study that tanker water suppliers indiscriminately extract water from certain groundwater sources, causing the water table to deplete beyond replenishable limits. Newspaper reports have also reported instances of informal water vendors often supplying untreated groundwater through tankers.[21, 22] It has been concluded in a survey report that due to the lack of adequate control to govern water quality and transparent delivery schedule, consumers pay much higher prices for water in the growing water tanker suppliers market.[23] Thus, to curb such unregulated activities and improve water supply through tanker trucks, the government is taking initiatives to develop a management system for tanker based water distribution in public-private partnership mode.[24] Such initiatives call for the proper integration of tanker water supply aspects and water resource management through improved planning and decision-making to ensure safe water access for all along with sustainable use of available water resources.

In this future direction, a design framework (long-term planning) is required to develop a route-map for improving current water supply situations by augmenting the pipe network supply with establishing tanker based water distribution systems in urban areas. This long-term planning is essential to support decisions for capital budgeting, treatment plant locations, treatment methods, number of tankers, tanker types etc. based on several factors such as (i) capital expenditure for tanker inventory and water treatment plants, (ii) available water source types in the city area, (iii) projections of water availability from those sources (iii) future demand assessments and demographic features of the urban area. Such decisions in long-term planning frameworks are based on the objectives such as maximizing return on investments, minimizing total annual cost etc., have a long-term impact on the cost of tanker water supply. On the other hand, a short-term planning framework is imperative to efficiently manage the day to day operations of an already established tanker based water supply systems for decisions concerning optimum utilization of treatment facilities, sustainable use of available water resources, ensuring water quality, timely delivery, preventing en-route water losses and demand-supply management through coordination among all groups in the supply chain to provide resilient water supply service in the current water-crisis scenario.[25] In this paper, we mainly focus on the short term planning aspects of an established water distribution network to ensure efficient and timely distribution of water to the consumers.

Towards this operational aspect in a short-term planning framework (although application dependent, typically ranges from hours- 1-2 weeks) research works describing optimal distribution planning for industrial gaseous products that are transported via trucks as bulk liquids are widely reported.[26-28] Furthermore, several research work present different approaches for time window discretization in the lateral domain of vehicles routing in supply-chain optimization problems.[29, 30] However, peculiarities involved in operating tanker based water supply system such as (i) maintenance of water quality while transportation, (ii) water availability from different sources, (iii) advanced water treatment requirements based on source water quality, (iv) rigorously incorporating time-dependent water consumption patterns of different consumers, (v) water demands for different purposes, e.g. drinking and domestic use, personal hygiene, medical use etc. (vi) tanker availability and suitability to various demand types etc. have not been addressed in the development of planning and scheduling framework at a commercial scale. Therefore, the optimization problem in this paper incorporates constraints based on a realistic abstraction of the above-mentioned complexities in the entire water supply chain (i.e., water source, treatment plants, and consumers) for the proposed short-term planning framework. In the purview of the short-term planning model herein, the water distribution problem is formulated as a MILP optimization framework to minimize the total operating cost for supplying treated water through a given tanker based water distribution system. The novel features and decision-making capability regarding several aspects of tanker water supply operations are demonstrated through the results of the proposed framework on a realistic urban case study in this paper. Furthermore, the presented framework can greatly simplify the computational complexity for solving real scale problems (which includes dozens of water intake points, several treatment facilities and thousands of consumers spanning across the city) to optimality with minimal relative optimality gap. The remainder of the paper is organized as follows: Section 2 provides a detailed description of complexities and challenges associated with a real scale planning problem of a tanker based water distribution system. Section 3 features the mathematical formulation of the proposed MILP optimization problem. The efficacy of the proposed framework is demonstrated through an illustrative realistic case study in Section 4 highlighting the optimal schedule of water tanker movements for water distribution, treatment etc. Section 5 then concludes with a summary of the novelties proposed in this paper and future research directions.

2. **Problem Description**

In this work, a tanker based water distribution system is considered to supply treated water from various water sources (e.g. freshwater source such as rivers, lakes etc. and groundwater source such as borewell, tube well etc.) to different type of consumers (e.g. households, commercial, hospitals, and institutions like schools etc.), given that both sources and consumers could span across the city. The consumers use water for multiple purposes, such as drinking, cooking, personal hygiene, washing, etc. Depending on the source as well as demand purposes, the water needs suitable treatment before supply (such as chlorine disinfection to prevent freshwater from microbial contamination, reverse osmosis for treating high TDS levels, removal of arsenic and other ions in groundwater, etc.).[31] Therefore, based on the water quality needs, this work considers three types of water states, viz. (i) untreated raw water from the source (ii) treated water for drinking purposes and (iii) treated water for domestic purposes; however, the formulation is flexible to include an arbitrary number of water states, representative of different water qualities. As shown in Figure 1, the entire water supply system can thus be considered as transportation of water products having states ($p \in P$) in tankers from a set of water sources ($s \in S$) to a set of water treatment plant ($s' \in WTP \subseteq S$), where water is treated up to standards of use for different purposes and then distributed among the set of consumers ($c \in C$). It is assumed here that water taken from the freshwater sources does not need any advanced treatment (e.g., reverse osmosis, UV/ozone disinfection) for all purposes other than drinking [19]. Thus, the tankers supplying water from freshwater sources can be considered adequate for all domestic usage with regular chlorine disinfection as an anti-microbial treatment. Subsequent sections of this paper describe the critical aspects related to treatment plant operation, distribution regions, and the required distribution capacity for efficient planning of tanker based water supply system.

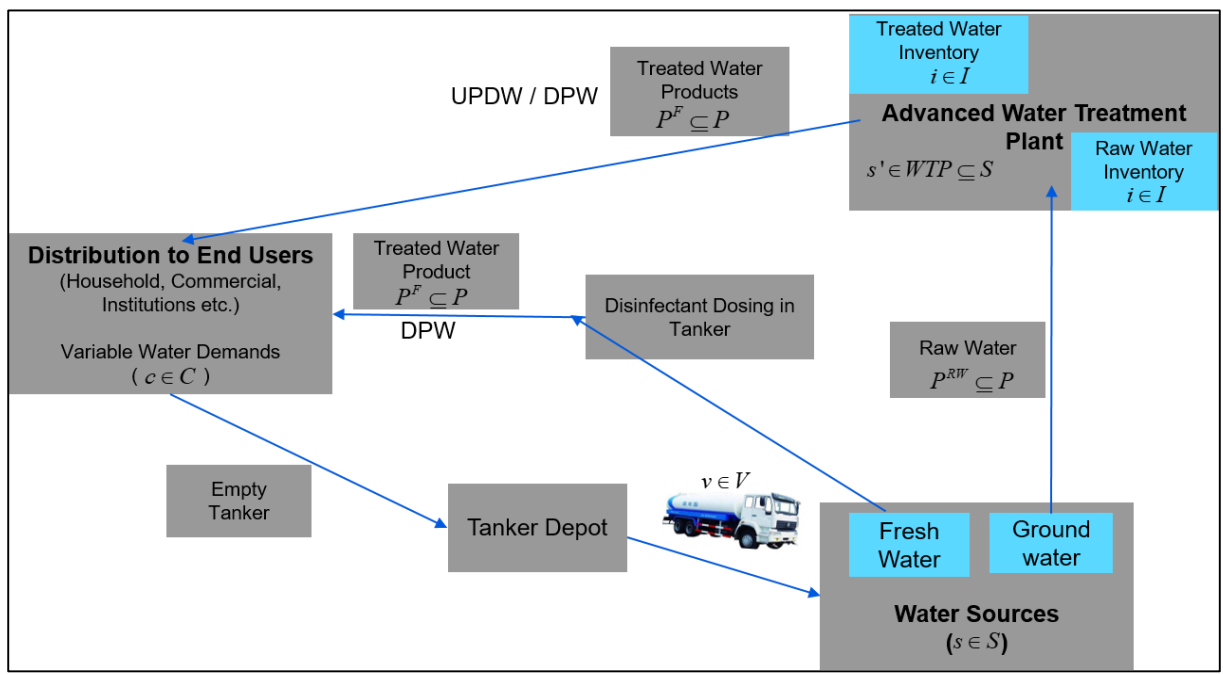

Figure 1: Design of Tanker based Water Distribution System

Considering the water treatments in treatment plants, separate inventories (*i*) are maintained for both untreated and treated water ($i \in I$). The water supplied to the treatment plant from the groundwater sources is constrained by the water extraction limits typically specified by the city planning board for groundwater reclamation.[32] In addition to this, the production of treated water in treatment facilities is also limited by factors such as processing throughput, reservoir capacity, and percentage recovery of treated water from the treatment process.[33] Furthermore, the operation at the treatment facilities is driven by the demands, so in the absence of demand, the treatment facilities can go through economical shutdown to minimize the operational costs. Therefore, if the treated water inventory at the treatment plant has sufficient water to fulfill the water demands, the operation of the treatment plant can be stopped to save energy costs.

To find an optimal schedule for the tanker routing, the total city area is proposed to be divided into several regions ($r \in R$) based on the water source, treatment plant, and consumer geographical locations. The reasoning behind such division of the total area into regions is to direct the optimizer for selecting tankers available near the sources. This will ensure effective utilization of the available tanker distribution capacity. Hence, the proposed framework takes input in the form of total number of tankers available in each region. Furthermore, if required, tankers can travel from one region to another depending on the suitability between sources-consumers, sources-water products and region-vehicle-water product. The source-consumer suitability parameter can be envisaged in the input data to the algorithm depending on the feasible distribution radius allowed to supply water from any source/treatment facility. Thus, this suitability matrix indicates a feasible mapping between the consumers and sources/ treatment facility. Similarly, source-water product suitability parameter can be envisaged based on the compatibility of the source/treatment facility to supply water up to the required purpose standards after appropriate disinfection or treatment respectively.

As shown in the Figure 1, there are four possible transportation sections for tanker movements in the entire water supply system: (i) transporting untreated water from groundwater sources to treatment plants, (ii) delivery of water from freshwater sources to consumers after proper disinfection in the tanker, (iii) delivery of treated water from treatment plant to consumer locations, and (iv) distribution of water to each consumer within that location. Therefore, depending on the transportation section, product *p* is supplied through a suitable type of tanker truck vehicles ($v \in V$) in the entire water supply chain. These tankers differ in their capacity and the materials used in the inner surface coatings. For instance, epoxy coated tanker trucks are required to supply clean and treated water to consumers, whereas usual stainless steel tanker trucks can be used for supplying untreated raw water to treatment plants [19]. On the other hand, some areas can be only catered by small capacity tanker trucks due to geographical constraints such as accessibility (narrow road infrastructure or elevated locations). Thus, based on the consumer-product, consumer-vehicle, and product-vehicle suitability, the selection of tanker trucks needs to be ensured.

Another important aspect of this tanker based water supply system is to rigorously capture the time-dependent water demands of different consumers in the resulting schedule. For example, commercial consumers require tanker water supply at time slots different than households. Thus, the water demands of different consumers need to be fulfilled at different periods in a day.[33] Further, catering to every individual consumer of a particular consumer type in a region requires detailed vehicle routing which will make the optimization problem intractable. Therefore, to make the problem solvable and computationally efficient, consumers of the same type are clustered together based on location within a region, with an assumption that once the vehicle enters the desired location; it can distribute water to each consumer within the given distribution time slot.

Furthermore, the transportation cost depends on several factors such as the selection of tanker type and capacity, distance to be traveled, transit time, etc. The travel time is calculated based on the distance to be covered and the average speed of the tanker truck (which will vary depending on tanker capacity). However, it is essential to note here that travel time alone is not sufficient to realistically represent the transit time for transportation between source to consumer, due to extra time required for other tasks such as (i) filling of water in the tanker, (ii) dosing the disinfectant and the time required for providing sufficient disinfection, and (iii) distribution time to each consumer. Hence, the actual transit time represented in this paper is calculated with due consideration to all of the above factors.

Considering the above-described complexities in a citywide tanker based water supply system, the proposed planning framework presented in this paper seeks to simultaneously optimize the water treatment and distribution aspects to enable the following decisions for transparent operation of the entire water supply chain:

(a) Generating the water treatment and distribution plan for at least one week
(b) Finding the optimal schedule for the movement of tanker trucks across the supply system (i.e., from the water source to treatment plants to consumers)

(c) Maximum utilization of the available distribution capacity to minimize distribution cost

(d) Exploration of integration between treatment plant operations, water availability at sources, and consumer demands to use the available water resources in a sustainable manner

(e) Maintaining the quality of water supplied to the consumers

(f) Deciding the optimal connection between the water source and consumer, which leads to minimum transportation costs considering the realistic transit time, water demand type and water demand timings

The detailed mathematical formulation of the above-mentioned optimization problem is presented in the following section.

*Remark 2.1*: The aspect related to water quality is addressed in the formulation from two perspectives: (i) fresh water sources are assumed to be fit for domestic/drinking water purposes after appropriate disinfection of water in tankers. In this perspective, we have assumed that given a proper disinfection time in our formulation, it can be assured to adhere to drinking water quality standards (whatever be the process used for disinfecting water in the tankers). Secondly, water from any ground water sources, being at a risk of contamination with metals/high TDS levels, it is first routed to treatment facilities for appropriate treatment in our formulation and then only supplied to consumers. In this perspective also, we have assumed that treatment method used at treatment facility is adequate enough to treat the water up to the required water quality standards for both DPW and UPDW purposes.

## 3. Mathematical Formulation

To develop a sufficiently rigorous model to assist in decision making at the planning layer and generate realistic targets for scheduling water tanker movements, a uniform and hourly discretization of the planning horizon is adopted in the proposed MILP optimization framework. The mathematical formulation of the problem described in the following encompasses four aspects: (i) water treatment plant operation, (ii) consumer demand fulfillment, (iii) tanker inventory utilization, and (iv) total operating cost, which forms the objective function. The notations of all indices, sets, parameters, and variables are provided in the nomenclature section along with the elucidation in the text as well. It is to be noted here that all the variables in the following MILP formulation are considered to be positive unless specified otherwise and colon symbol (:) is used to represent "such that", as typically used in mathematical formulations.

### 3.1 Water Treatment Plant Operation

As mentioned in the preceding section, the proposed formulation considers the option of stopping the water treatment plant operation when there is no demand which needs to be fulfilled by treating groundwater. Therefore, a binary variable '$yOp_{s,t}$' is considered which assumes the value of '1' in the time period when the treatment plant is operating to produce the required product type, else it is set to be '0'. Also, there are constraints during the operation of the treatment plant such that it needs to continue the operation for a minimum uptime period (denoted by parameter $T_s^{UT}$). Similarly, the treatment plant requires specific down time period (denoted by parameter $T_s^{DT}$) before restarting the operation. These operational transition constraints for treatment facility are sought to be incorporated by considering positive continuous variables '$xSUp_{s,t}$' and '$xSDn_{s,t}$' in combination with binary variable '$yOp_{s,t}$' as shown in Equations (1-4), curtailing the total number of binary variables in the formulation. As shown in Figure 1, three types of water sources are considered in this paper, 1) freshwater (FW), 2) groundwater (GW) and 3) treatment facilities (TF) and the parameter $STy_s$ is used to identify these source types. Hence, the constraint using parameter $STy_s = $ 'TF' are only applicable to treatment facility type sources. The parameter $Op_{s,t}^{ini}$ in Equation 2 represents the initial operational state of the treatment facility at the start of the current planning horizon.

$$xSUp_{s,t} - xSDn_{s,t} = yOp_{s,t} - yOp_{s,t-1}$$
$$\forall s(:STy_s = 'TF'), t(:t = 2..NT) \qquad (1)$$

$$xSUp_{s,t} - xSDn_{s,t} = yOp_{s,t} - Op_{s,t-1}^{ini}$$
$$\forall s(:STy_s = 'TF'), t(:t = 1) \qquad (2)$$

$$\sum_{k=t-T_s^{UT}+1:(t-T_s^{UT}+1)>0}^{t} xSUp_{s,k} \leq yOp_{s,t}$$
$$\forall s(:STy_s = 'TF'), t(:t=1..NT) \quad (3)$$

$$\sum_{k=t-T_s^{DT}+1:(t-T_s^{DT}+1)>0}^{t} xSDn_{s,k} \leq 1 - yOp_{s,t}$$
$$\forall s(:STy_s = 'TF'), t(:t=1..NT) \quad (4)$$

Furthermore, different consumers require water of different qualities for different purposes which demands various treatment standards, e.g., groundwater needs treatment up to 500 ppm TDS for drinking purposes, while for domestic purposes, water is suitable if treated only up to 1000 ppm TDS [31]. In the proposed framework, we are considering water with different quality/purposes as separate states. In this framework all the water states are denoted using index $p \in P$. Accordingly, two subsets ($P^{RW} \subseteq P$ and $P^F \subseteq P$) are considered respectively for raw water state and treated water states/final products which are of deliverable quality. However, the treatment plants can treat raw water for only one final water product state at a time. This constraint related to the selection of water product to be produced by a treatment plant at any time period is captured in Equation (5). The binary variable '$yPSl_{s,p,t}$' is used to denote the production of a specific final product state $p$ at a time period $t$ in treatment facility $s$. Here the parameter $SP_{s,p}$ takes the value 1 if source is suitable to supply product $p$.

$$\sum_{p(:SP_{s,p}=1 \& p \in P^F)}^{P} yPSl_{s,p,t} = yOp_{s,t}$$
$$\forall s(:STy_s = 'TF'), t(:t=1..NT) \quad (5)$$

Now, mass balance constraints for treatment facility are included to account the consumption of raw water received at treatment plant and production of final treated water products. Consequently, the raw water inventory balance at the treatment plant is written as shown in Equation (6-7), where the variable '$xQ_{s,i,p,t}$' captures the quantity of water state $p$ available in inventory $i$ of source $s$ at time $t$ and the variable '$xSSupl_{s',s,p,t}$' indicates quantity of water state $p$ that is supplied from the water source to the treatment plant. Here the parameter $T^{RWTransit}_{s',s,p}$ denotes time taken in filling/emptying tanker and travel from source to treatment facility. Furthermore, the parameter $STpt_s$ in Equation 6 denotes throughput of the water treatment plant (also known as treatment capacity) in kiloliter per hour. It is to be noted here that any treatment operation produces treated water products at a fixed percentage recovery[34], which is stored as a fraction in the parameter '$\beta_{s,p}$' used in the following Equations 6-7. Therefore, LHS of equation 6 equates the amount of water present in the inventory at time $t$ ('$xQ_{s,i,p,t}$') with the terms on RHS which can be interpreted as the quantity of water present in the inventory at time $t-1$ ('$xQ_{s,i,p,t-1}$') and the amount of water added to the inventory at time $t$ (which is supplied from the source at $t-T^{RWTransit}_{s',s,p}$ time steps ago, taking $T^{RWTransit}_{s',s,p}$ travel time from ground water source to treatment plant) and subtracting the quantity of water used from the raw water inventory for producing treated water as term 3. The parameter '$Q^{ini}_{s,i,p}$' in Equation 7 represents the amount of the raw water available in the inventory at the start of the planning horizon

$$xQ_{s,i,p,t} = xQ_{s,i,p,t-1} + \sum_{s'(:SS_{s',s}=1\,\text{and}\,t-T^{RWTransit}_{s',s,p}>0)}^{S} xSSupl_{s',s,p,t-T^{RWTransit}_{s',s,p}}$$
$$- \sum_{p'(:p'\in P^F)}^{P} yPsl_{s,p',t} * STpt_s * \frac{1}{\beta_{s,p'}}$$
$$\forall\, s(:STy_s = 'TF'), i(:i = 'RWI'),$$
$$p(:p \in P^{RW}\,\text{and}\,SIP_{s,i,p} = 1), t(:t = 2..NT) \tag{6}$$

$$xQ_{s,i,p,t} = Q^{ini}_{s,i,p}$$
$$\forall\, s(:STy_s = 'TF'), i(:i = 'RWI'),$$
$$p(:p \in P^{RW}\,\text{and}\,SIP_{s,i,p} = 1), t(:t = 1) \tag{7}$$

Further, the raw water gets treated in the treatment plant and stored in the treated water inventory (TWI). The treated water from inventory is then delivered to multiple consumers. Therefore, a variable '$xDeCon_{s,c,p,t}$' is used to capture the contribution of a source $s$ in fulfilling the demand of the water product of state $p$ for any consumer $c$ and accordingly mass balance around the treated water product inventory in the treatment plant is described in Equation 8-9. Thus, LHS of equation 8 represents quantity of water present in treated water reservoir at time $t$ and equates this in RHS with three terms denoting quantity of treated water present in the treated water reservoir at time ($t$-1) as '$xQ_{s,i,p,t-1}$', quantity of water treated and added to the reservoir at time $t$ ($STpt_s$) and subtracting the total water supplied to consumers from this reservoir at time $t$ ($xDeCon_{s,c,p,t}$).

$$xQ_{s,i,p,t} = xQ_{s,i,p,t-1} + yPSl_{s,p,t} * STpt_s - \sum_{c}^{C} xDeCon_{s,c,p,t}$$
$$\forall s(:STy_s = 'TF'), i(:i = 'TWI'),$$
$$p(:p \in P^F \text{ and } SIP_{s,i,p} = 1), t(:t = 2..NT) \quad (8)$$

$$xQ_{s,i,p,t} = Q_{s,i,p}^{ini}$$
$$\forall s(:STy_s = 'TF'), i(:i = 'TWI'),$$
$$p(:p \in P^F \text{ and } SIP_{s,i,p} = 1), t(:t = 1) \quad (9)$$

Nextly, both the raw water and treated water inventory at the treatment plants are constrained to be always maintained under the physical bounds of minimum ($ICap_{s,i,p}^{min}$) and maximum storage ($ICap_{s,i,p}^{max}$) capacities of the reservoir. These constraints are given by Equations (10-11).

$$xQ_{s,i,p,t} \geq ICap_{s,i,p}^{min}$$
$$\forall s(:STy_s = 'TF'), i,$$
$$p(:p \in P \text{ and } SIP_{s,i,p} = 1), t(:t = 1..NT) \quad (10)$$

$$xQ_{s,i,p,t} \leq ICap_{s,i,p}^{max}$$
$$\forall s(:STy_s = 'TF'), i,$$
$$p(:p \in P \text{ and } SIP_{s,i,p} = 1), t(:t = 1..NT) \quad (11)$$

Along with the satisfaction of minimum and maximum limits on inventory, treatment plant operation also requires maintenance of buffer capacity in raw water inventory to allow smooth water pumping for treatment (i.e. the quantity of water that is needed to be always maintained in the reservoir as below which the water cannot be pumped easily from the reservoir to the treatment section). Similarly, treated water inventory are also usually bound to maintain target volumes for meeting emergency water demands, upcoming maintenance periods etc. Therefore, these buffer and target capacity constrints on treatment plant water inventories are also included in this framework as described in Equations 12-13 respectively. In these equations, the variables '$xBCV_{s,i,p,t}$' and '$xTV_{s,i,p,t}$' denote the amount of violation from buffer capacity limit and target respectively. The superscripts (+) and (-) on variable '$xTV_{s,i,p,t}$' in Equation 13 represents positive and negative violations respectively from target limit of a treated water inventory. It is to be noted here that violation from buffer capacity and any negative violation from target capacity to satisfy these constraints is subsequently penalized by appropriate factors in the total operating cost.

$$xQ_{s,i,p,t} \geq ICap^{buffer}_{s,i,p} - xBCV_{s,i,p,t}$$
$$\forall s(: STy_s = 'TF'), i(: i = 'RWI'),$$
$$p(: p \in P^{RW} \text{ and } SIP_{s,i,p} = 1), t(: t = 1..NT) \quad (12)$$

$$xQ_{s,i,p,t} - xTV^{+}_{s,i,p,t} + xTV^{-}_{s,i,p,t} = ICap^{Target}_{s,i,p,t}$$
$$\forall s(: STy_s = 'TF'), i(: i = 'TWI'),$$
$$p(: p \in P^{F} \text{ and } SIP_{s,i,p} = 1), t(: ICap^{Target}_{s,i,p,t} > 0) \quad (13)$$

where, the parameters '$ICap^{buffer}_{s,i,p}$' denotes the required buffer capacity in raw water inventory and $ICap^{Target}_{s,i,p,t}$ indicates inventory target to be achieved at time period $t$.

Lastly, the aspect related to raw water extraction limit from any groundwater source at a time period $t$ is formulated as Equation 14 in which the parameter '$SMax_s$' stores the maximum permissible water extraction limit from a groundwater source $s$. This regulatory limit by Ground Water Control Board prevents exploitation of the ground water sources beyond replenish-able limits. This limit is usually expressed as KL/hr. Thus, '$SMax_s$' is not a time sensitive parameter.

$$xSSupl_{s,s',p,t} \leq SMax_s$$
$$\forall s(: SMax_s > 0), s'(: SS_{s,s'} = 1),$$
$$p(: p \in P^{RW}), t(: t = 1..NT) \quad (14)$$

**3.2 Consumer Demand Fulfillment**

As described in the previous section, consumers of the same type in a location (household/commercial/institutional etc.) are clustered in one group and based on the water consumption patterns of each consumer type, the daily water demands of all consumers in a group are aggregated product wise. The maximum and minimum demand of a customer cluster for a product at time period $t$ is captured by parameters '$De_{c,p,t}^{\max}$' and '$De_{c,p,t}^{\min}$' respectively. Thus, equation 15 represents the consumer demand fulfillment constraint describing that the total water supplied from all the sources to a consumer in a time period $t$ should be greater than the minimum consumer demand at that time period. Equation 16 represents similar constraint imposed due to parameter $De_{c,p,t}^{\max}$. Also, it can be seen as indicated in equations 15-16 that multiple sources can contribute (captured by variable '$xDeCon_{s,c,p,t}$') in fulfilling the demand from a consumer cluster, wherein the suitability of a source to cater to a certain consumer cluster is captured by the parameter $SC_{s,c}$ and a similar source-product suitability condition is indicated by the parameter $SP_{s,p}$. Here the parameter $T_{s,c,p}^{Transit}$ denotes the time taken for tanker filling, disinfection process (if it is fresh water source) and travel from source to consumer. It is also important to note here that restricting equations 15-16 for time instants having $De_{c,p,t}^{\min} > 0$ and $De_{c,p,t}^{\max} > 0$ implies that there is no need of considering such constraint on the variable $xDeCon_{s,c,p,t}$ for time periods where $De_{c,p,t}^{\min} = 0$ and $De_{c,p,t}^{\max} = 0$ respectively.

$$\sum_{\substack{s(:SP_{s,p}=1 \text{ and } SC_{s,c}=1 \\ \text{and } t-T_{s,c,p}^{Transit}>0)}}^{S} xDeCon_{s,c,p,t-T_{s,c,p}^{Transit}} \geq De_{c,p,t}^{\min}$$

$$\forall c,\ p(: p \in P^F),\ t(: De_{c,p,t}^{\min} > 0) \qquad (15)$$

$$\sum_{\substack{s(:SP_{s,p}=1 \text{ and } SC_{s,c}=1 \\ \text{and } t-T_{s,c,p}^{Transit}>0)}}^{S} xDeCon_{s,c,p,t-T_{s,c,p}^{Transit}} \leq De_{c,p,t}^{\max}$$

$$\forall c, p(: p \in P^F), t(: De_{c,p,t}^{\max} > 0) \qquad (16)$$

### 3.3 Tanker inventory

Typically, different types of tankers (in terms of capacity and inner surface material compatibility) are required to supply different water products across the water supply chain. Therefore, the optimization in the planning layer should also account for effective utilization of the available capacity of the vehicles to transport water products in different sections, respecting the delivery timings and transit time required to reach the destination locations. As briefly mentioned in section 2, four types of movements are possible for the considered water states in the tanker water supply system: (i) raw water supply from groundwater source to treatment plants (ii) treated water state/final product distribution from treatment plants (secondary source) to consumer locations (iii) treated water distribution from freshwater sources to consumer locations after chlorine-based disinfection in tanker and (iv) distribution of water to each consumer in the consumer groups once the tanker reached the consumer location. However, instead of detailed distribution routing of the tanker for each consumer in the consumer group, the planning layer optimization in this work considers a sufficient distribution time period for the fourth section while calculating total time capacity in subsequent equations. This assumption along with the concept of overall time capacity balance as further explained in equation 19 makes the optimization problem solvable and computationally efficient without losing on the accuracy of timely delivery. The constraints on the utilization of available tankers capacity to fulfill total demands are described next.

In this direction, first term on RHS of equation 17 denotes total quantity of water that is to supplied from a source to consumers in the entire time horizon. Therefore, this quantity should be equal to the total available tanker capacity, which is represented by second term on RHS of equation 17. Thus, the variable '$xCDistb_{s,c,p,v}$' is used in equation 17 to calculate the amount of treated water product $p$ that should be distributed from source $s$ using vehicle type $v$ suitable to both product $p$ and consumer $c$ to fulfill its total demand. Extending the same concept to the quantity of water that has to be supplied from a source to treatment facility in the total time horizon, equation 18 is written with a variable '$xVSSupl_{s,s',p,v}$' is to calculate the amount of raw water that should be transported in vehicle type $v$ from a source $s$ to a treatment facility (secondary source) $s'$. In summary equation 17-18 represents the constraint that total distribution need should be equal to the total available distribution capacity.

$$\sum_{t=1}^{NT} xDeCon_{s,c,p,t} - \sum_{v(:CPV_{c,p,v}=1)}^{V} xCDistb_{s,c,p,v} = 0$$

$$\forall c, p(: p \in P^F), s(: SC_{s,c} = 1 \& SP_{s,p} = 1) \quad (17)$$

$$\sum_{t=1}^{NT} xSSupl_{s,s',p,t} - \sum_{v(:SSPV_{s,s'p,v}=1)}^{V} xVSSupl_{s,s',p,v} = 0$$

$$\forall s, s'(: SS_{s,s'} = 1), p(: p \in P^{RW}) \quad (18)$$

Furthermore, based on the average speed and distance between the source and consumer, the parameter $T^{Travel}_{s,c,p,v}$ denotes only the travel time taken by vehicle type $v$ to travel from source $s$ for delivering product $p$ to consumer $c$. Similarly, $T^{Travel}_{s,s',p,v}$ stores only the travel time taken by vehicle type $v$ to travel from source $s$ for supplying product $p$ to treatment plant $s'$. However, in addition to the travel time, total transit time for distributing a product $p$ from source $s$ to consumer $c$ also includes (i) time required to prepare tanker for water supply which includes tanker travel from depot to source, tanker filling etc. (ii) time required in disinfection process of fresh water to prevent microbial re-growth in filled tanker water and (iii) time required for distributing water to each consumer in the group, once the tanker reached the consumer group location. With these terms represented by the parameters $T^{Prep}_{s,v}$, $T^{Disf}_{s,v}$ and $T^{Distb}_{c}$ respectively, Equation 19 describes overall time capacity balance. Here, the term "time capacity" denotes the time required for supplying unit quantity of water product. Therefore, as written in first term in LHS of equation 19, multiplication of time capacity with total quantity of water to be distributed with the available tankers capacity ($xCDistb_{s,c,p,v}$) will give the total time required to supply water from source to consumers. Similarly, second term represents the total time required to supply water from sources to treatment facilities. And, now this total required time to supply water products should be less than or equal to total time at disposition in the planning horizon with the available tankers and any extra hired tankers in the region $r$. This is represented by first term on RHS of equation 19. It is to be noted here that the time involved in calculating time capacity in first term is the total time needed to travel to the desired location, distribute the required quantity of water state $p$ from source $s$ using vehicles type $v$ and return to the originating depot. Similarly, for supplying water from source to treatment facilities, second term includes preparation, travel and return time. Here, the parameter, $VA_{r,v,p}$ denotes the number of vehicles of type $v$ available in region $r$ that can transport water of state $p$. Furthermore, in case the capacity of the available tanker inventory is insufficient to perform all the transportation needs in a region, extra vehicles can be purchased/hired incurring premium costs. The continuous variable '$xVExQ_{r,v,p}$' in Equation 19 denotes the extra tanker capacities that are required to be purchased/hired, and the corresponding costs are minimized as a part of the objective function. Also, important to note that based on the suitability parameter '$RVP_{r,v,p}$' tankers can move from one region to any other suitable region $r$ for transporting water product of state $p$ and return back to originating region.

$$\sum_{c(:SP_{s,p}=1 \text{ and } RS_{r,s}=1 \text{ and } CPV_{c,p,v}=1)}^{C} \left( \frac{2*T_{s,c,p,v}^{Travel} + T_{s,v}^{Prep} + T_{s,v}^{Disf} + T_{c}^{Distb}}{VQ_v} \right) * xCDistb_{s,c,p,v}$$
$$+ \sum_{s(:RS_{r,s}=1 \text{ and } SSPV_{s,s',p,v}=1)}^{S} \left( \frac{2*T_{s,c,p,v}^{Travel} + T_{s,v}^{Prep}}{VQ_v} \right) * xSTf_{s,s',p,v}$$
$$\leq NT * \left( VA_{r,v,p} + \frac{xVExQ_{r,v,p}}{VQ_v} \right)$$
$$\forall r,v,p(:RVP_{r,v,p}=1) \quad (19)$$

Although Equation 19 appropriately represents constraint related to the supply of all products $p$ in a region $r$ using the total capacity of vehicle types $v$ in the total available time period, it might provide an overestimation of the transportation capacity available in a region at a particular time period. This can be understood as the concept of overall time capacity balance will give correct estimates of transportation capacity available in a region on a particular time period, only if demands are uniformly sparse in the total time horizon. However, if demands are concentrated in some particular hours the available transportation capacity will be overly estimated as equation accounts the balance of tanker capacity for total time period. Hence, in addition to overall time capacity balance in Equation 19, we further impose an hourly transportation capacity balance for tankers in equation 24. Towards this, the calculation of tanker capacity consumed in a time period either for supplying water from source/treatment plant to consumers or from source to treatment plants is required. This calculation is shown in Equation 20-21 for the delivery of products from source/treatment plant to consumers. The variable '$xPDl_{s,c,p,v,t}$' is introduced here to denote the quantity of treated water product $p$ delivered using vehicle type $v$ from source $s$ to the consumer $c$ in a time period $t$. Thus, the term on RHS of equation 20 represents total quantity of water product to be supplied from source to consumers in a time period $t$. And, this total quantity should be equal to the total available transportation capacity in that particular time period. This transportation capacity is represented by the LHS of equation 20.

$$\sum_{v(:CPV_{c,p,v}=1)}^{V} xPDl_{s,c,p,v,t} = xDeCon_{s,c,p,t}$$

$$\forall s, c(:SC_{s,c}=1), p(:SP_{s,p}=1 \text{ and } p \in P^F), t(:t=1..NT) \tag{20}$$

On the other hand, the term on LHS of equation 21 represents quantity of water to be supplied from source to consumers in the total time horizon. Therefore, this should be equal to total available transportation capacity of type $v$ which is suitable to consumer $c$ and product $p$, as represented by RHS of the equation 21.

$$\sum_{t=1}^{NT} xPDl_{s,c,p,v,t} = xCDistb_{s,c,p,v}$$

$$\forall s, c(:SC_{s,c}=1), p(:SP_{s,p}=1 \text{ and } p \in P^F), v(:CPV_{c,p,v}=1) \tag{21}$$

Similarly, Equation 22-23 describes the calculation of tanker capacity used in supplying raw water from the source to treatment plants in a particular time period. Here, the variable '$xRw_{s,s',p,v,t}$' indicates the amount of raw water state $p$ supplied from source $s$ to treatment plant ($s'$) using vehicle type $v$ in time period $t$.

$$\sum_{v(:SSPV_{s,s',p,v}=1)}^{V} xRw_{s,s',p,v,t} = xSSupl_{s,s',p,t}$$

$$\forall s, s'(:SS_{s,s'}=1), p(:SSPV_{s,s',p,v}=1), t(:t=1..NT) \tag{22}$$

$$\sum_{t=1}^{NT} xRw_{s,s',p,v,t} = xVSSupl_{s,s',p,v}$$

$$\forall s, s'(:SS_{s,s'}=1), p(:SSPV_{s,s',p,v}=1), v \tag{23}$$

Thus, Equation 20-23 calculates the amount of water that needs to be transported using a vehicle type $v$ at time period $t$. Now, using the Equation 24, the total transportation capacity ('$xVQ_{r,v,p,t}$') of tanker type $v$ that is available at any time period $t$ for transporting water of state $p$ in region $r$ is calculated, as represented by first term on LHS. This is then equated with terms in RHS for hourly balance as following: (i) first term of RHS denotes transportation capacity available at time period $t$-1, (ii) second term subtracts transportation capacity consumed in supplying water from sources to consumers in time period $t$ (iii) third term subtracts transportation capacity consumed in supplying water from source to treatment facilities in time period $t$, (iv) fourth term makes addition in the transportation capacity due to tankers returning back to depot in same time period $t$ after supplying water to consumers and (v) fifth term also makes addition in the capacity due to tankers returning back to the depot in same time period $t$ after supplying raw water at treatment facilities. Hence, equation 24 appropriately balances available tanker capacity limitations in each time period. This equation also ensures that the empty tankers after distribution return to the tanker depot in the originating region, where it gets added to the available transportation capacity for the next trip in the next time period.

$$xVQ_{r,v,p,t} = xVQ_{r,v,p,t-1}$$
$$- \sum_{c}^{C} \sum_{s(:SP_{s,p}=1 \& RS_{r,s}=1 \& SC_{s,c}=1)}^{S} xPDl_{s,c,p,v,t}$$
$$- \sum_{s(:RS_{r,s}=1)}^{S} \sum_{s'(:SSPV_{s,s',p,v}=1)}^{S} xRw_{s,s',p,v,t}$$
$$+ \sum_{c}^{C} \sum_{s(:SP_{s,p}=1 \& RS_{r,s}=1 \& SC_{s,c}=1)}^{S} \sum_{\substack{t-1(:TE_{t-1}>(TE_{t'}+\frac{2*T^{Travel}_{s,s',p,v}+T^{Prep}_{s,v}+T^{Disf}_{s,v}+T^{Distb}_{c}}{24})) \\ t'(:TS_{t-1}\leq(TS_{t'}+\frac{2*T^{Travel}_{s,s',p,v}+T^{Prep}_{s,v}+T^{Disf}_{s,v}+T^{Distb}_{c}}{24}))}} xPDl_{s,c,p,v,t'}$$
$$+ \sum_{s(:RS_{r,s}=1)}^{S} \sum_{s'(:SSPV_{s,s',p,v}=1)}^{S} \sum_{\substack{t-1(:TE_{t-1}>(TE_{t'}+\frac{2*T^{Travel}_{s,s',p,v}+T^{Prep}_{s,v}}{24})) \\ t'(:TS_{t-1}\leq(TS_{t'}+\frac{2*T^{Travel}_{s,s',p,v}+T^{Prep}_{s,v}}{24}))}} xRw_{s,s',p,v,t'}$$
$$\forall r,v,p,t(:t=2..NT)$$

(24)

where, the initial capacity balance of tanker inventory is given by Equation 25. Accordingly, it is important to note here that, the presented framework aids in making decisions of extra tanker capacity requirements ("$xVExQ_{r,v,p}$") for meeting the total demands in complete planning horizon at the start of the planning horizon itself.

$$xVQ_{r,v,p}^{ini} = VQ_v * VA_{r,v,p} + xVExQ_{r,v,p}$$
$$\forall r, v, p(: RVP_{r,v,p} = 1) \quad (25)$$

### 3.4 Objective Function

The main aim of this formulation is to minimize the operating cost for supplying water through tanker based system to provide timely delivery of water to consumers while utilizing the available water resources, treatment facility and tanker inventory optimally. Hence, the objective function (*xobj*) is formulated as shown in Equation 26. The first two terms in the RHS of Eq.26 signifies the total transportation cost of tankers for distributing treated water from sources to consumers and supplying raw water from sources to treatment facilities respectively. Further, the third and fourth term denotes penalty cost on violating target and buffer capacities in treatment plants respectively, while the last term represents penalty cost on extra tanker vehicles purchased/hired.

As mentioned in the previous section that the tanker travel cost parameter is calculated as cost per quantity based on the mileage and capacity of each tanker type. Therefore, multiplication of travel cost parameter "$TrCost_{s,c,p,v}^{Distb}$" with total quantity of water to be supplied from sources to consumers in first term on RHS of equation 26 will give the total cost incurred in supplying water in this section of supply system. Similarly, second term represents total cost incurred in supplying water from source to treatment plants. However, the penalty cost for hiring extra tankers is usually given on the basis of number of tankers. Therefore, the last term calculates the number of tankers required to fulfill extra capacity needs and then multiplies by the cost factor.

$$xobj = \sum_{s}^{S}\sum_{c}^{C}\sum_{p(:SP_{s,p}=1)}^{P}\sum_{v}^{V} Tr\text{Cost}_{s,c,p,v}^{Distb} * xCDistb_{s,c,p,v}$$

$$+ \sum_{s}^{S}\sum_{s'(:SS_{s,s'}=1)}^{S}\sum_{p}^{P}\sum_{v(:SSPV_{s,s',p,v}=1)}^{V} Tr\text{Cost}_{s,s',p,v}^{RWsupply} * xVSSupl_{s,s',p,v}$$

$$+ \sum_{s}^{S}\sum_{i}^{I}\sum_{p}^{P} TV\text{Cost}_{s,i,p} * \sum_{t=1}^{NT} xTV_{s,i,p,t}$$

$$+ \sum_{s}^{S}\sum_{i}^{I}\sum_{p}^{P} BCV\text{Cost}_{s,i,p} * \sum_{t=1}^{NT} xBCV_{s,i,p,t}$$

$$+ \sum_{r}^{R}\sum_{v}^{V}\sum_{p(:RVP_{r,v,p}=1)}^{P} VEx\text{Cost}_{v,p} * \frac{xVExQ_{r,v,p}}{VQ_{v}}$$

(26)

## 4. Results and Discussion

This section discusses the efficacy of the proposed MILP optimization framework towards handling the peculiarities of planning and management of tanker based water supply system, through an application on a representative case study. This case study is based on the practical insights about typical tanker water supply system operations from *Just Paani Water Supply Solutions.*[35] Such enterprises that operate tanker water supply typically manage a network of tanker water vendors in different areas of the city. These vendors supply water either from ground water or fresh water sources in different capacity tankers. The consumers and vendors are mapped with each other by the enterprise company based on (i) the customer demand (quantity and time slot) data, (ii) delivery locations and (iii) availability of water tankers with several vendors, in the network across the city. However, in the present situation, due to the lack of any planning and scheduling framework for real-time operation of the tanker water supply system, the suppliers face challenges of untimely water delivery and water quality monitoring issues. The model input data for the application case study presented in the following section is qualitatively reconciled with insights from *Just Paani Water Supply Solutions*, literature survey studies and newspaper articles/ reports on tanker water supply market to represent real world relevance and scale of such planning and scheduling problem.[11, 20, 23]

The results of framework are demonstrated to assist decision makers on various aspects of tanker water supply system operation to ensure quality water and fulfillment of all costumers demands, in the required turnaround times. Also, a typical scenario of a maintenance period for the treatment facility is designed to present the ability of proposed formulation in assisting towards optimal decisions. The developed formulation is programmed and solved in FICO Xpress Optimization[36] environment using "mmxprs" module version 2.8.1 on 16 GB RAM machine with Intel i7 (3.6 GHz) processor. Depending on the urban settlement, there can be several treatment plants and water sources catering to the city water supply. Therefore, the formulation presented above is programmed generically and is flexible for application to other problem sizes for a satisfactory solution in terms of the achieved optimality gap. The details of the case study design and scenarios are described as following:

**4.1 Application Case Study**

This case study considers two fresh water sources (FW2, FW3), three groundwater sources (GW1, GW3, GW4) and two water treatment plants (TF1, TF2) to serve the water demands of the city as shown in Figure 2. Based on the location of these sources and consumers, the entire city is divided into three regions (R1, R2, and R3) having a different combination of water sources and treatment plants as depicted by Figure2. It is to be noted here that (i) such kind of regional division (North/east/west/south) is typical for water supply aspects in many Indian cities[37] and (ii) in order to include a wide spectrum of all possible physical scenarios, we have also included the aspect of one zone (R3) having multiple water sources but no treatment facilities yet established. Furthermore, consumers in each region are clustered in three types, namely, household consumers (HHC), commercial consumers (CC), and hospital and health care institutions (HC). Amongst these consumers, household and commercial consumers are assumed to have water demands for both domestic purposes (DPW) and ultra-pure drinking water (UPDW), whereas hospital consumers water demand is assumed for UPDW type only. The raw water from all the groundwater sources is first transported to either of the treatment facilities in region R1 and R2 and then treated water is delivered to consumers in different locations. To transport water in different sections, this case study considers tankers in two capacities 6000 L (6T) and 10,000 L (10T) for raw water whereas, treated water is supplied using tankers of three capacities (3000 L, 6000

L and 10,000 L). Also, considering that the typical geographical locations of household communities are narrow areas, it is assumed that water is delivered to HHC by only 3T and 6T tankers while the other two consumer types can be supplied by all three tanker sizes. A total number of available tankers of each type and size is quoted region wise in the algorithm input file as provided in supplementary material with this paper.

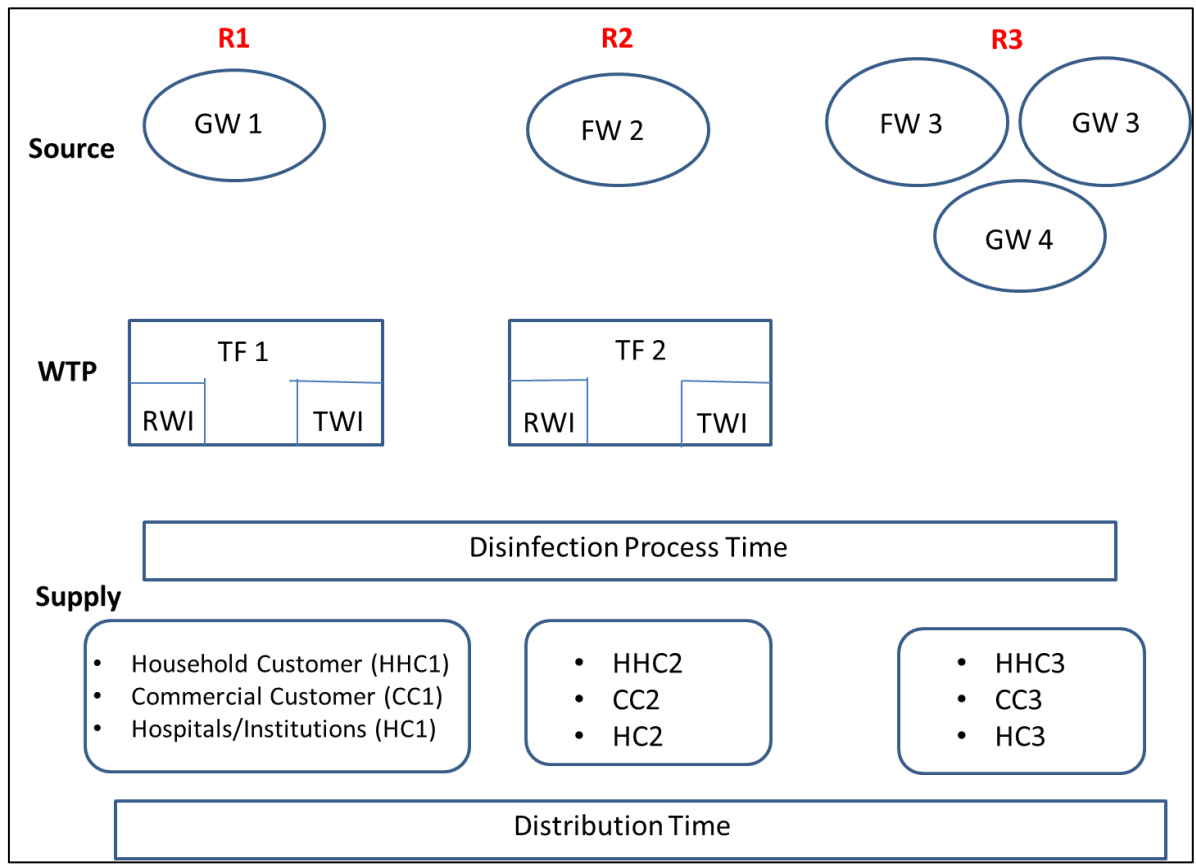

**Figure 2: Description of Components in the Application Case Study**

The water demands data in this case study is designed to fulfill the typical tanker water demand of 11.5 MLD in three regions of city, as shown in Table 1. As discussed briefly in earlier sections that the water demands of each consumer type have two aspects – quantity and timings. The water demands of each consumer type are different in quantity as well as at different time periods of the day. Considering both these aspects, as shown in Table 2 and Table 3 we have divided the total water demands of each consumer type into two delivery slots (morning & evening) according to the general urban water consumption pattern.[13, 33] It is to be noted from Table 2 and Table 3 that there is a significant difference in water demand pattern on weekdays and weekends for household and commercial consumer types.

The planning horizon in this case study is assumed to be of 8 days, demnd horizon of 7 days and remaining details of algorithm input data for this case study can be referred from supplementary information file.

Table 1: Assumed water demands of each consumer type in three regions

| Consumer | Water Demand in R1 (KL) | Water Demand in R2 (KL) | Water Demand in R3 (KL) |
|---|---|---|---|
| HHC | 330 | 480 | 600 |
| CC | 2600 | 3250 | 3250 |
| HC | 375 | 250 | 300 |

Table 2: Total Water Usage Pattern of each Consumer Type

| Consumer | Morning | | Evening | |
|---|---|---|---|---|
| | Weekdays | Weekends | Weekdays | Weekends |
| HHC | 60% | 80% | 40% | 20% |
| CC | 70% | 45% | 30% | 55% |
| HC | 70% | 50% | 30% | 50% |

Table 3: Water Demand Timings and Distribution Slot for each Consumer Type

| Consumer | Morning | | Evening | |
|---|---|---|---|---|
| | Weekdays | Weekends | Weekdays | Weekends |
| HHC | 5 AM - 9 AM | 6 AM - 10 AM | 4 PM - 8 PM | 2 PM - 6 PM |
| CC | 9 AM – 12 PM | 10 AM – 1 PM | 5 PM – 8 PM | 4 PM – 7 PM |
| HC | 6 AM - 8 AM | 6 AM - 8 AM | 6 PM - 8 PM | 4 PM - 6 PM |

Table 4: Optimization problem statistics

| Problem Characteristics | Value |
|---|---|
| Constraints | 22493 |
| Continuous Variables | 42771 |
| Binary Variables | 384 |
| Time (s) | 43.9 |
| Optimality gap (%) | 0.06 |

The problem statistics for this case study are shown in Table 4. The following results demonstrate various characteristics of the framework which are helpful for planning and management of available water resources, tanker inventory and treatment plant operators to provide reliable service of clean and safe water to the public using tanker based water supply system. The values reported in all the figures and results are normalized to keep the scale of graphs clear and concise. However, actual data corresponding to one normalized unit in each figure is provided along with this discussion.

Taking the values of the variable '$xDeCon_{s,c,p,t}$' for each source, Figure 3(a) is drawn to show the percentage contribution of each source in fulfilling the total demand (DPW and UPDW) of all consumers in a weekday. As seen in Figure 3(b), this percentage distribution is different on weekend days due to changes in water demand patterns. Thus, the framework can optimize source utilization based on the water demands.

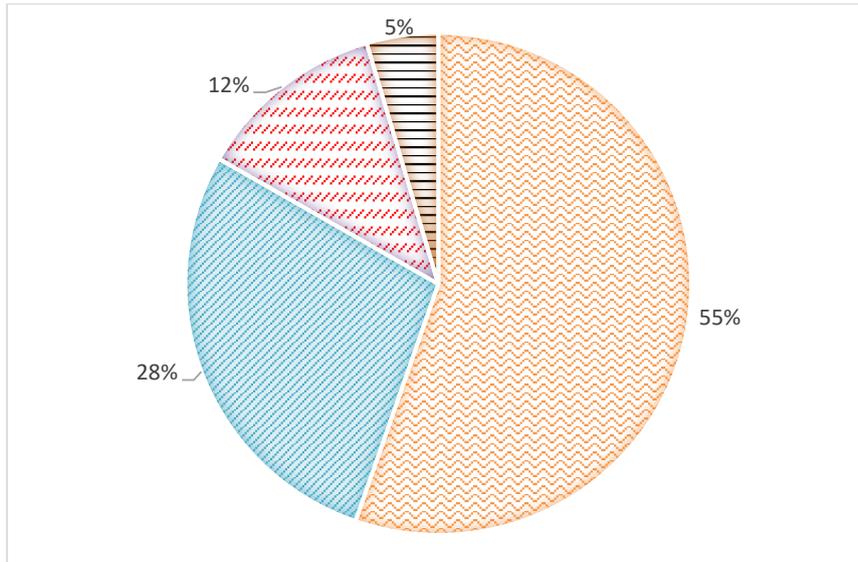

(a)

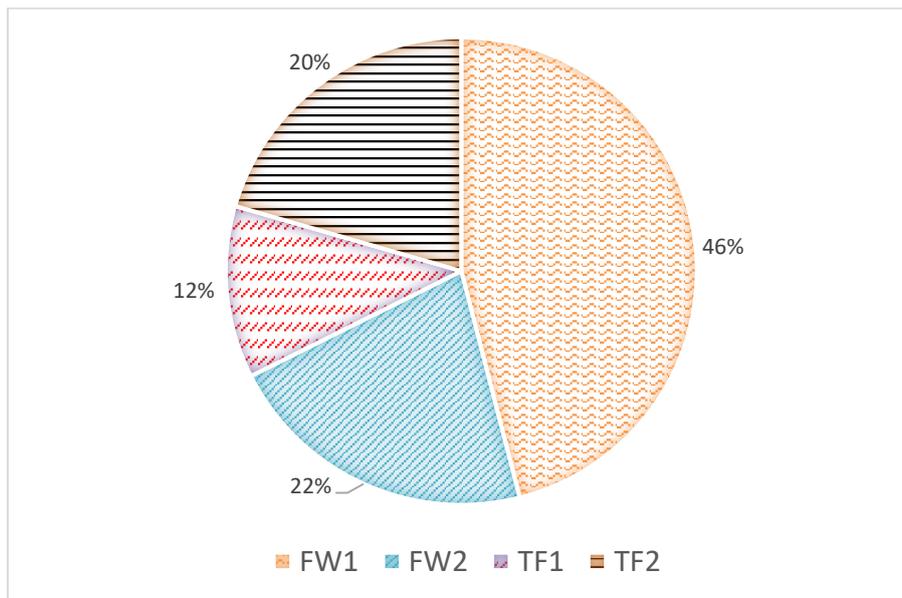

(b)

**Figure 3: Contribution of each Source in total Demand Fulfillment of a day (a) Weekday (b) Weekend**

As mentioned earlier that the water treatment facilities maintain separate reservoirs for raw water received from groundwater sources and treated water for supply to consumers. Thus, corresponding to the day 4 operation, raw water received at TF1 in each hour from several ground water sources and further getting treated at a constant throughput rate (given by $STpt_s$) is shown in Figure 4. The blue, orange and green bars in the Figure 4 represents amount of water supplied from GW1, GW3 and GW4 respectively and sum of these quantity constitutes total water received at TF1 in time period *t*. Furthermore, the yellow bars show the quantity of water used from raw water inventory to treat upto DPW quality and thus black line with y-axis on RHS describes the change in water levels in the raw water inventory at the end of each hour. This is known as raw water inventory profile. Also, Figure 4 shows the necessity of planning frameworks due to complexity of problem in hand where a treatment plant in one region is receiving raw water from groundwater sources of multiple regions in different time periods to fulfill total demands while honoring all the operational constraints. It is to be noted here that some of the following figures are shown for one day or one consumer type to depict the features of the framework while keeping the description of results concise and clear. However, similar observations can be drawn from the results file for each consumer and complete planning horizon. The normalization value of raw water quantity corresponding to one unit on y-axes in Figure 4 is 430 kilo Liters.

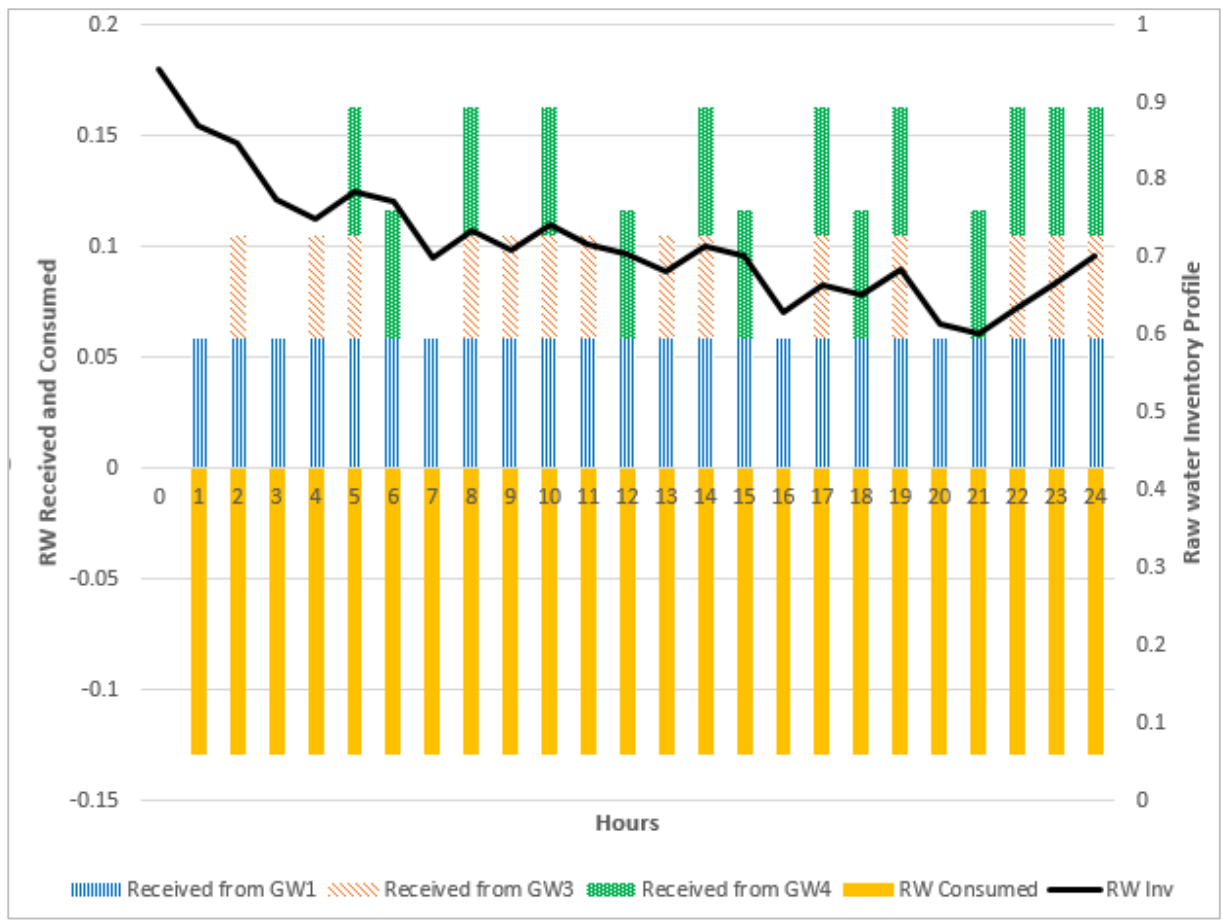

**Figure 4: Hourly Details of Raw Water Inventory in TF1 for day 4 operation**

Similarly, Figure 5 illustrates the daily treated water inventory profile for DPW production in TF1. It can be seen from inventory profile in Figure 5 that the WTP operation starts from day 1, while water demands are starting from day 2 for consecutive 7 days in the total 8 days planning horizon of this case study. Furthermore, it can be observed from the Figure 5 that although the quantity of water to be supplied from TF1 is same on all weekdays, the framework is optimizing the daily production of DPW while respecting the minimum, maximum and target bounds on the inventory profile. The normalization value of DPW water quantity corresponding to one unit on y-axes in Figure 5 is 850 kilo Liters.

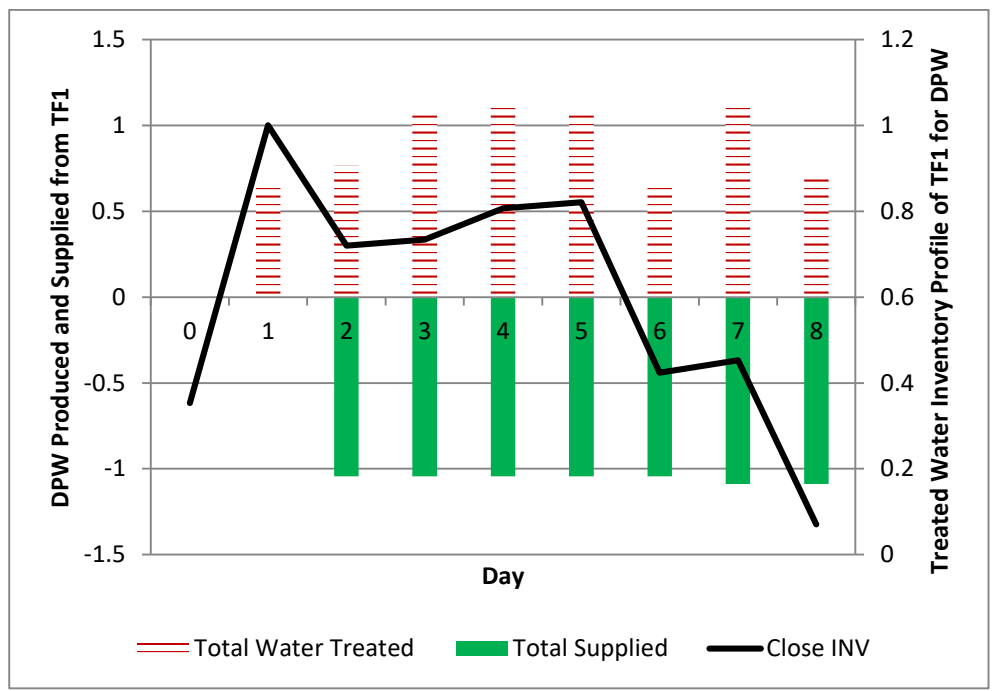

**Figure 5: Daily Production and Supply from Treated Water (DPW) Inventory**

Another feature of the framework to cater the demands of multiple consumers in different regions from one source is shown in Figure 6. Thus, Figures 5 and 6 combinely shows the simultaneous optimization of treatment plant operation and timely delivery aspects for the smooth operation of all components in the entire water supply system. The normalization values of amount of water delivered on day 3 from each water source in Figure 6 (FW3, FW2, TF1, TF2) are (6.325 MLD, 3.201 MLD, 1.387 MLD, 0.433 MLD) respectively.

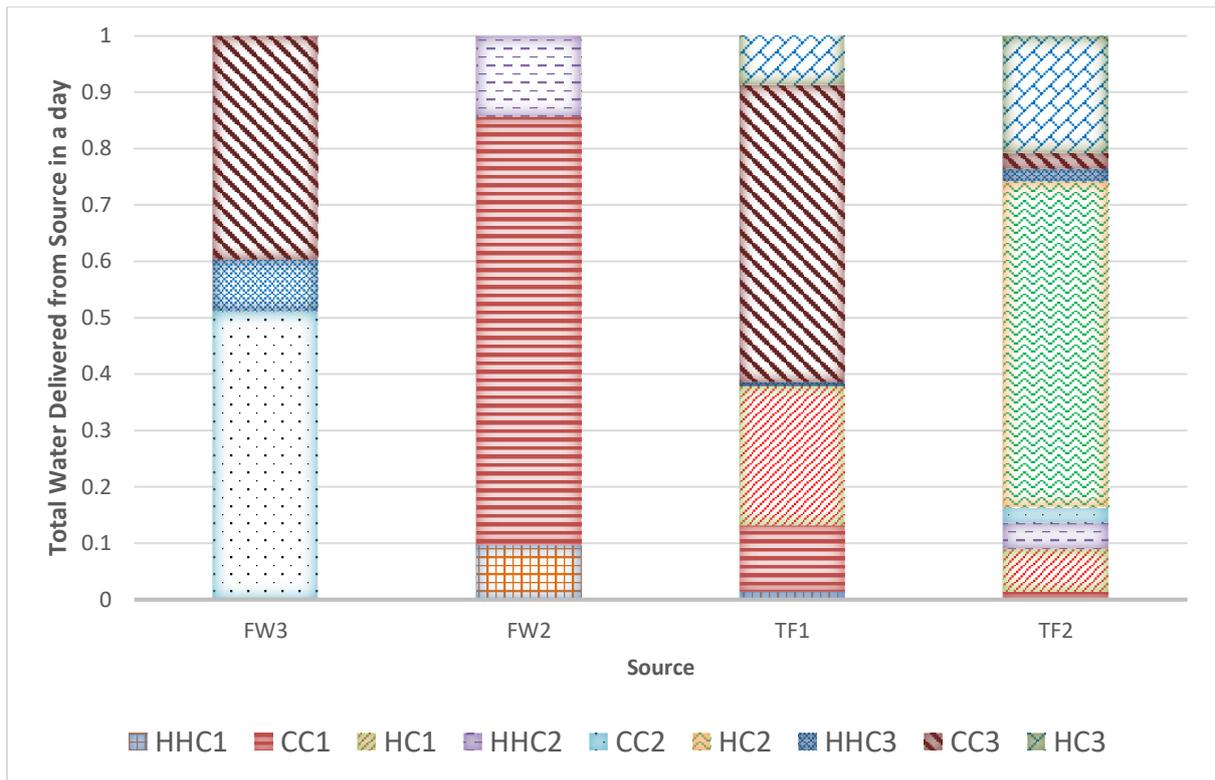

Figure 6: Each Water Source Catering to Multiple Consumers in different regions

Further to understand the water supply service to each consumer type, Figure 7 shows the demand fulfillment pattern of DPW product for consumer CC1, in both morning and evening time periods. This figure shows that in situations of increased water demands which cannot be fulfilled by the sources present in that region, the supply can be planned to be met by sources of another region in an optimal manner. For example, as seen from the demand pattern in Table 2, there is high water demand in evening periods compared to morning periods for CC1. To meet this extra demand, TF1 is also supplying water to CC1 on day 7 and day 8 evening periods. Thus, Figure 7 highlights the framework's feature where a consumer is serviced by multiple sources depending on the quantity of water demand, demand timings, and tanker availability in the region. The normalization value of CC1 consumer DPW water demand in Figure 7 is 2.594 MLD.

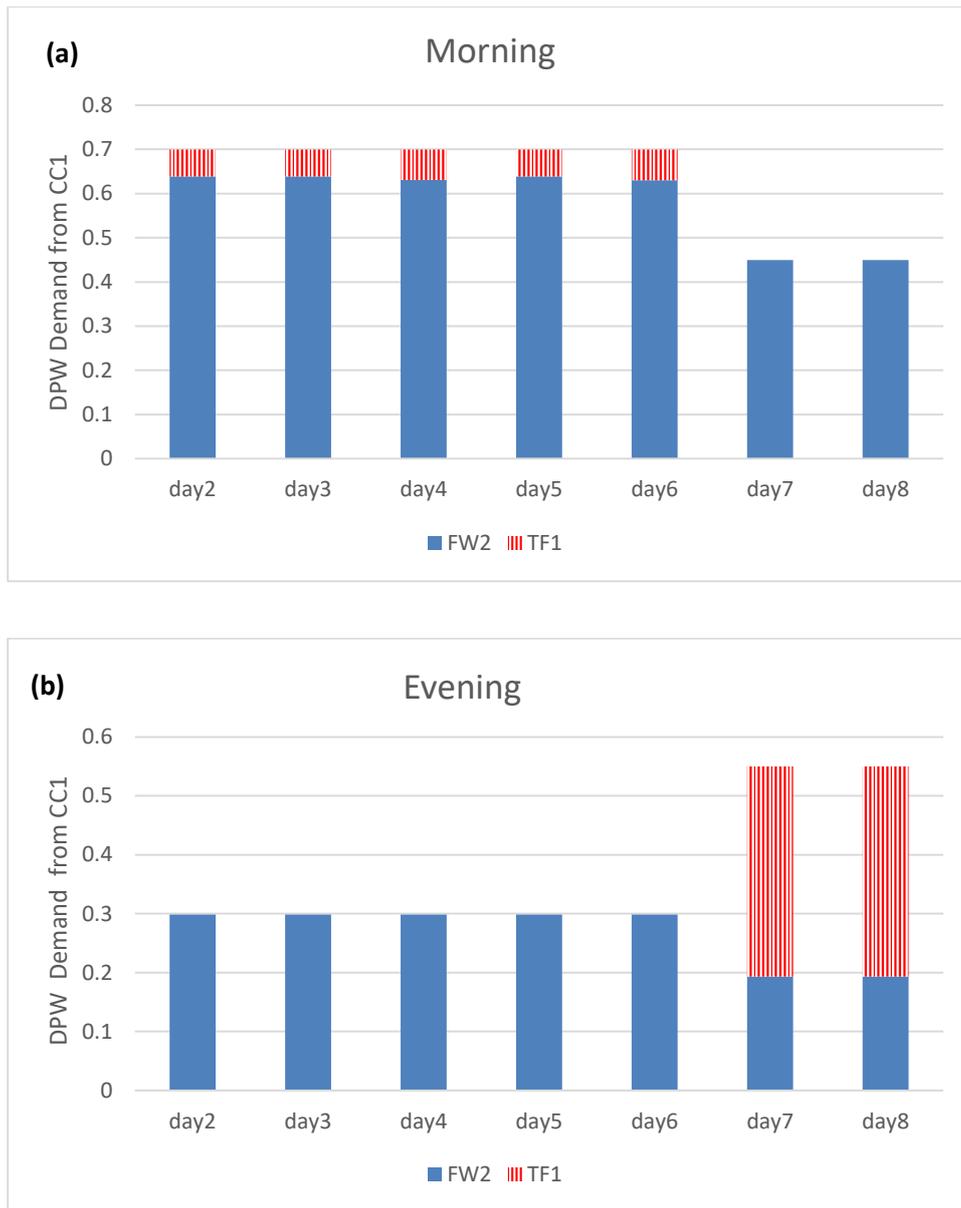

Figure 7: Daily Demand fulfillment pattern of CC1 Consumer (a) Morning period (b) Evening period

As discussed earlier that the three different capacity tankers (3T, 6T, and 10T) are considered in this case study. The maximum use of available capacity respecting the constraints of consumer-tanker suitability, transit time, and delivery timings are other vital aspects of optimization in the proposed planning framework. Figure 8 depicts that according to the water demand quantity in different time periods, tankers of various capacities are selected to supply water at minimum transportation cost. Thus, the framework provides a detailed understanding of questions such as how much quantity of water to be delivered from different suitable sources to a consumer by optimally utilizing multiple tanker capacities. The normalization value of CC2 DPW morning period demand in figure 8 is 2269 kilo liters.

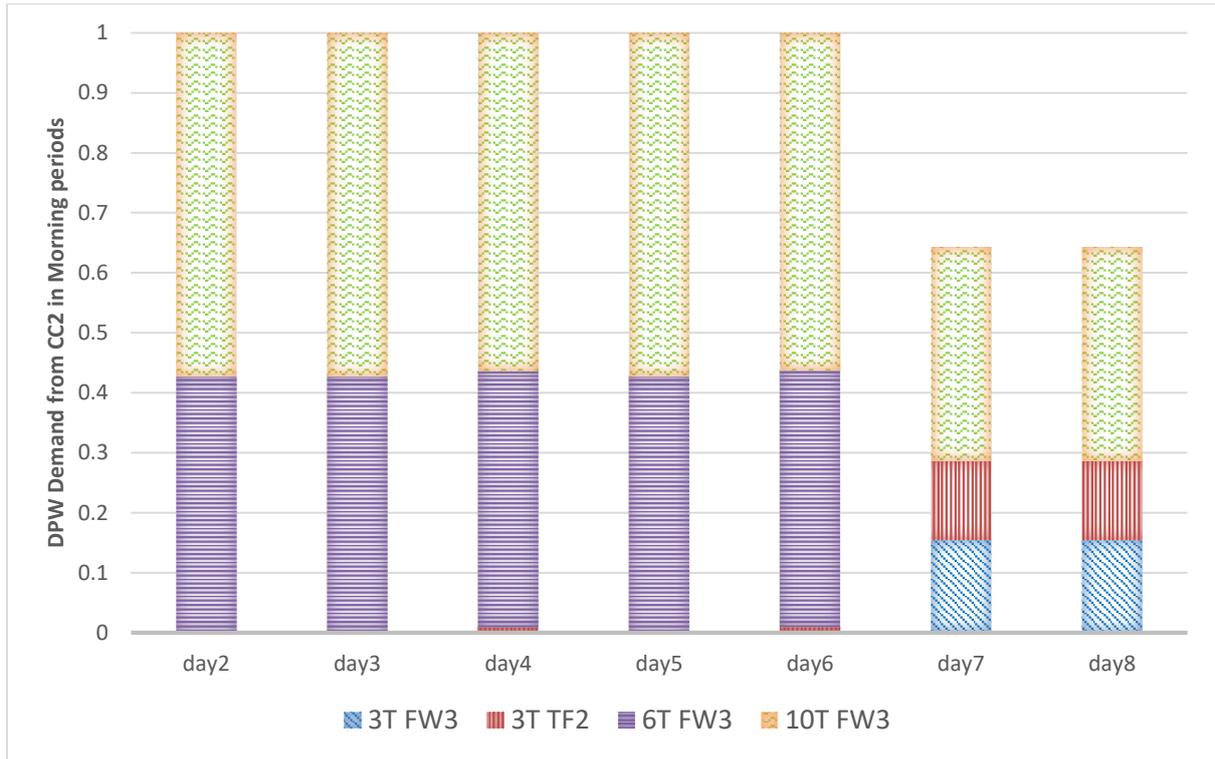

**Figure 8: Source-Consumer-Truck association for daily demand fulfillment of CC2 consumer type**

Table 5 shows the distribution details of initially available tankers, extra tankers required to be purchased/ hired to fulfill total demand and number of trips taken by each tanker in the 8-day horizon period. The interpretation of fraction values of total number of trips in Table 5 should be drawn as for e.g. if a total of 267 tankers of 3T type (having capacity of 3000L) are making 8.34 trips, then 176 tankers [267-(0.34*267) =176] are making 8 round trips while 91 tankers (=0.34*267) are completing 9 trips in the total horizon. It can be further observed from Table 5 that the larger capacity trucks are taking more trips than smaller capacity. This is also intuitive as transportation cost per unit quantity of water supplied by 6T and 10T tanker is much lower as compared to the 3T tanker. This in turn indicates the potential to improve by minimizing total annual capital cost in a long-term planning framework in future studies.

**Table 5: Distribution Details of Tanker Availability and total trips in the planning horizon**

| Truck Type | Product | Truck Capacity (L) | Total Truck Availability | Extra Truck Purchased/Hired | No. Of Trips |
|---|---|---|---|---|---|
| 3T | DPW | 3000 | 260 | 7 | 8.34 |
| 3T | UPDW | 3000 | 60 | 0 | 4.96 |
| 6T | DPW | 6000 | 515 | 163 | 10.10 |
| 6T | UPDW | 6000 | 80 | 0 | 5.66 |
| 10T | DPW | 10,000 | 235 | 0 | 10.7 |
| 10T | UPDW | 10,000 | 40 | 0 | 8.75 |
| 10T | RW | 10,000 | 190 | 0 | 9.52 |

Lastly, Figure 9 shows the percentage distribution of various cost components in the objective function. It can be observed from Figure 9 that corresponding to the optimal solution, penalty cost for buffer and target limit violations in treatment plant inventories is negligible compared to the tanker travel costs for supplying water from sources to consumers and from sources to treatment facilities.

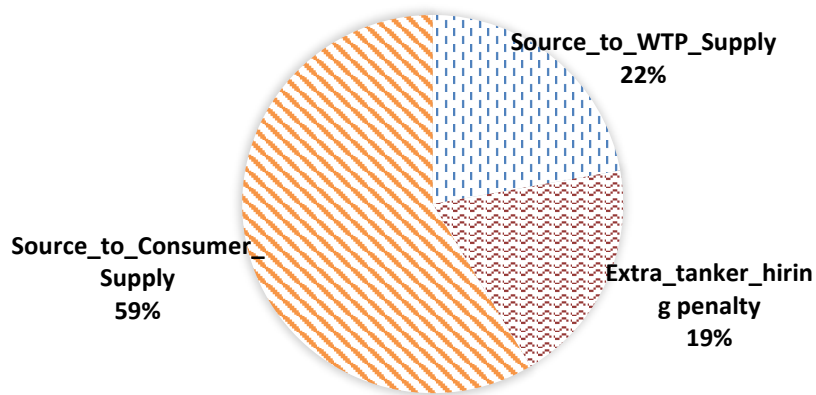

**Figure 9: Percentage distribution of cost components in objective function corresponding to optimal solution**

## 4.2 Water Treatment Plant Maintenance Period during the planning horizon

This scenario describes the situation when the water treatment plant is under maintenance duration for several practices such as reservoir cleaning, pumps maintenance, servicing of treatment installations, etc. In such instances, the other region treatment facilities need to be operational to provide resilient water supply to consumers of those regions. This feature is shown here by shutting down the operation of 'TF1' for the first two days of the planning horizon (1-48 hours). The demand pattern of consumers in all regions remains the same as considered in the previous section case study. The resulting optimal schedule of TF1 and TF2 operation (generated by the variable '$yOp_{s,t}$') is shown in Figure 10(a) and compared with the optimal schedule of normal scenario (Figure 10 (b)) described in the previous section. This comparison shows that such water inventory and demand constraints during such planned maintenance periods can be efficiently taken care of by optimizing treatment operations for complete planning horizon in both TF1 and TF2.

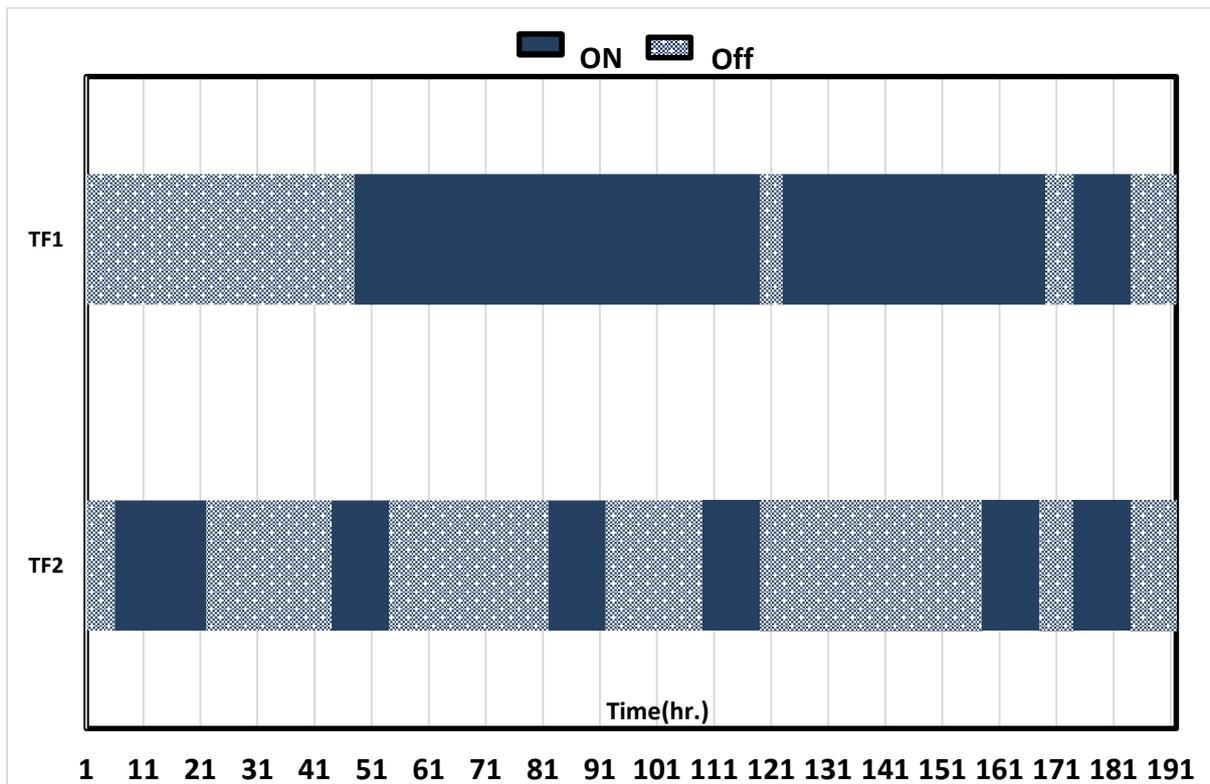

(a)

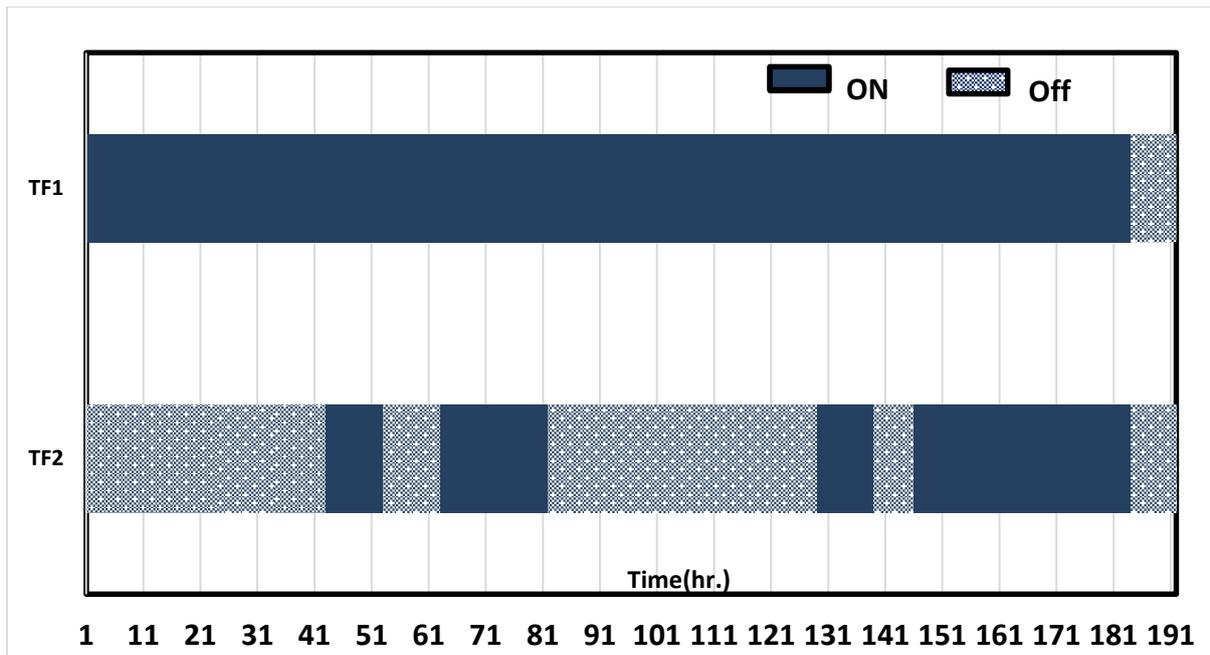

(b)

**Figure 10: Optimal Schedule of Treatment Facilities operation (a) Maintenance Period (b) Normal case**

Further to understand the water supply aspects and quantity of water to be treated in treatment facilities for consumer demand fulfillments during such maintenance period, TWI profile of UPDW in TF2 and TF1 is compared in this scenario in Figure 11 (a) and (b) respectively. The comparison shows the framework's feature to provide optimal water treatment targets for continuing resilient water supply during such scenarios, such that operating cost for the entire planning horizon can be minimized while respecting constraints related to increased demand load on TF2 during TF1 maintenance, inventory capacities, and timely delivery. The normalization value corresponding to one unit on y-axis in Figure 10 (a) is 1140 kilo Liters and 565 kilo Liters in Figure 10 (b). Also, Table 6 is shown to indicate the overall increase in total objective function value in such maintenance period compared to the normal case. The significant increase in penalty cost on hiring extra tankers in maintenance period is expected as in such maintenance scenario where Region 1 has no other treatment plant to treat water from its ground water sources, model is seeking to fulfill region 1 consumer demands from a farther distance at a minimal cost. Therefore, while 10T type tanker trips and cost of supplying water from GW1 to TF1 is reduced, more number of 6T type tankers are hired for this period and demands are being fulfilled from another region. These details of type of extra tankers purchased as well as change in total number of trips are also compared for both cases in Table S11 (supplementary material).

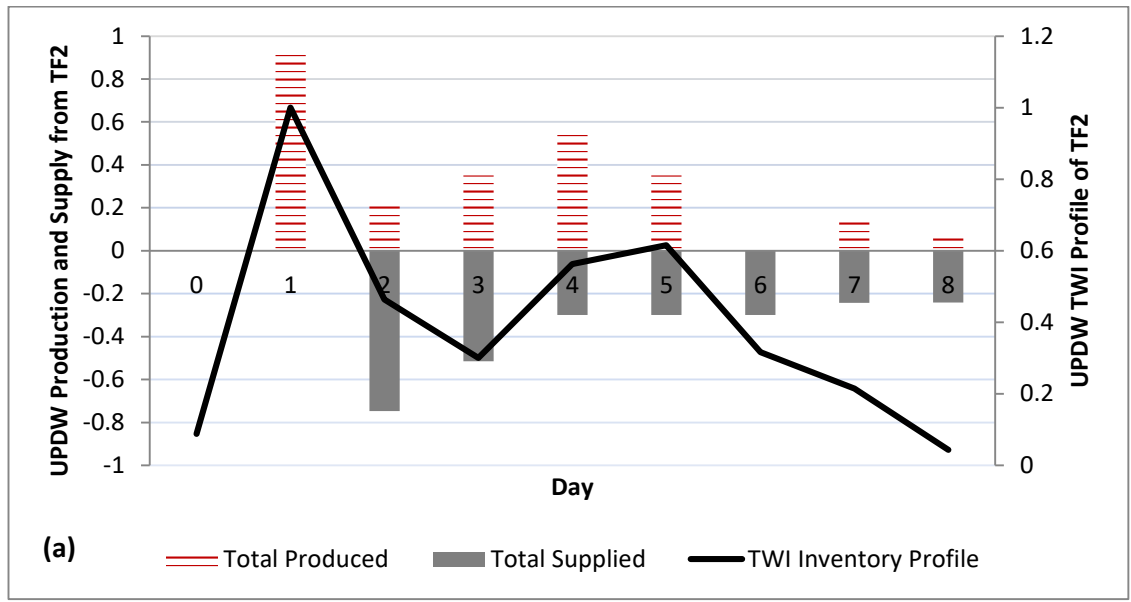

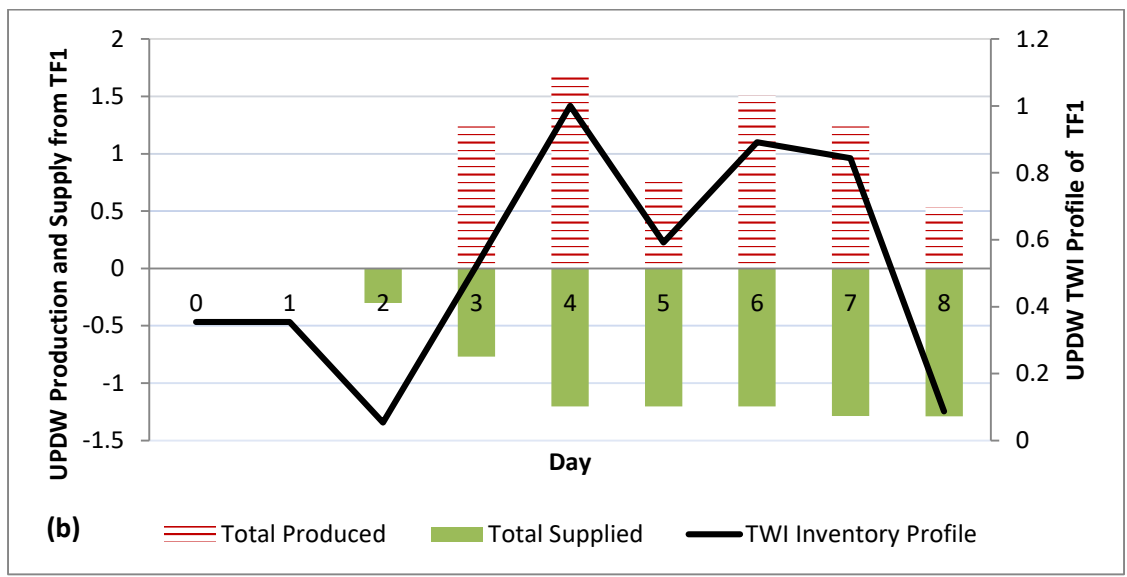

**Figure 11: Comparison of Treated Water Inventory Profile of UPDW during the Maintenance Period in planning horizon: (a) TF2 (b) TF1**

Table 6: Objective function and cost comparison for normal and maintenance periods

| Scenario Cost (INR) | Cost of Source to Consumer Supply | Cost of Source to WTP supply | Penalty Cost of buffer capacity violations | Penalty Cost of Target limit violations | Penalty Cost of extra tankers hire/purchase | Total Operating cost ($x_{obj}$) |
|---|---|---|---|---|---|---|
| Normal | 10041300 | 3798600 | 26141.7 | 32050 | 3300000 | 17198091.7 |
| TF1 Maintenance | 10045600 | 2604500 | 27297.2 | 32050 | 6258330 | 18967777.2 |

*Remark 4.1*: As observed through results in this section, the presented framework deals not only water tanker transportation scheduling but decisions on several other aspects such as treatment plant operations, selection of water treatment type for different purpose water demands, extra tankers hiring to fulfill water demands etc. Therefore, it is important to note here that there are no two different time scales of planning these operations and scheduling tanker movements. Rather, the tanker scheduling decisions are integrated in the planning horizon itself. Hence, this is termed as a short-term planning framework, where the thrust of the study is to explore integration of treatment facility operations, water demands and distribution logistics at a common time scale for efficient management of tanker water supply system. This is different than the typical industrial supply-chain optimization problems where production planning is explicitly carried out on a much larger time horizon compared to distribution scheduling.

*Remark 4.2*: The application case study assumes planning for water demand of 11.5 MLD (Million liters per day) in all 7 days of demand horizon. While daily variations in this demand value might present in a realistic scenario, the assumed value is considered as a nominal demand data for the presented deterministic framework. However, water usage patterns are different on weekdays and weekends as given in Table 2 of the manuscript. Thus, in totality the aspect of time varying demand is present in the data. Furthermore, the data used in the presented case study is designed to be just a simple representative of all possible scenarios, while formulation in section 3 is quite generic to work with any kind of real-data.

## 5. Conclusion

This paper presents a short term planning and scheduling framework for augmenting access to clean and treated water through tanker based water distribution system in urban areas. The main objective of the proposed framework is to find the optimal tanker movement schedule (from sources to treatment plants to consumers) while minimizing the total operating cost for timely delivery of water to consumers and sustainable use of available water resources. The proposed MILP optimization formulation rigorously incorporates major constraints related to operational nuances of different components (viz. water sources, treatment facilities, and consumers) in tanker based water supply system. Therefore, the developed planning framework respects various limitations of water treatment, transit time, and distribution aspects to provide the optimal target schedule, which is achievable in representative applications. The framework proposed in this paper is applied on case studies, which incorporate the intricacies of various real-world scenarios. The results presented in this paper conclusively demonstrate several features of the proposed framework, essential to manage the operation of tanker based water supply system such as (i) contribution of each water source to fulfill total demand (ii) consumers catered with respect to each water source (iii) treatment plants facilitating water treatment from groundwater sources in different regions (iv) reservoir profiles of raw water and treated water inventory at treatment plants (v) operation schedule of treatment plant (vi) association of tanker types, consumer and sources in optimal tanker delivery schedule (vii) decisions for hiring extra tankers at the start of operation to fulfill total water demands. Furthermore, quantitative analysis in the case study shows that the travel cost of tankers for delivering water from sources to consumers constitiutes major component, however penalty cost of hiring extra tankers is the most sensitive component in the total operating cost. Future work has the scope to extend this study for developing a long term planning framework and integrating with this short-term model for operating water supply through tankers in urban areas. Also, this formulation can be extended for tanker supply system in rural areas which are not connected to water mains in the pipe network, as well as peri-urban regions which suffer problems of head-loss in the pipe distribution network.

**Supporting Information:**

The supplementary material associated with this manuscript contains (i) data used for designing application case study, (ii) Gant chart showing specimen of tanker movements schedule in one day of planning horizon and (iii) example of expanding constraint equations in the formulation. This information is available free of charge via the Internet at http://pubs.acs.org/.

**Author's Information**

Corresponding Author
*E-mail: ravigudi@iitb.ac.in.
ORCID:
 Ravindra D. Gudi: 0000-0002-1684-4174

**Acknowledgment**

The authors gratefully acknowledge funding to the first author from the Ministry of Human Resources Development, Government of India. We also extend special thanks to Mr. Liju James from *Just Paani Water Supply Solutions, Bangalore, India,* for providing realistic insights of typical tanker water supply system case study presented in this paper.

**Acronyms:**

DPW   Domestic Purpose water

FW    Fresh water Source

GW    Ground water Source

RW    untreated raw water from ground water source

RWI   raw water inventory at treatment facility

TF   Treatment facility

TW   treated water in treatment facility

TWI   Treated water inventory at treatment plant

UPDW   Ultra-pure Drinking water

**Nomenclature**

The following symbols are used in this paper:

**Sets:**

| | |
|---|---|
| $C$ | set of consumers |
| $I$ | set of water inventory at treatment plants |
| $P$ | set of water product states |
| $P^{RW} \subseteq P$ | set of raw water states |
| $P^{F} \subseteq P$ | set of treated water state/final products |
| $R$ | set of regions for water supply in urban area |
| $S$ | set of water sources |
| $V$ | set of tanker vehicles |
| $WTP \subseteq S$ | set of water treatment plants |

**Indices:**

| | |
|---|---|
| $c \in C$ | consumer |
| $i \in I$ | water inventory in the treatment facility |
| $p \in P$ | water product state |
| $r \in R$ | region |
| $s, s' \in S$ | water source |
| $v \in V$ | tanker truck vehicle type |
| $t$ | time period |

**Parameters**

$BCV\text{Cost}_{s,i,p}$ penalty cost for violating the buffer capacity of product state $p$ in inventory $i$ at source $s$

$CPV_{c,p,v}$ assumes value equal to 1 if vehicle $v$ is compatible to deliver product of state $p$ to consumer $c$

$De_{c,p,t}^{\min}$ minimum aggregate demand from consumer group $c$ at time $t$ for water product of state $p$

$De_{c,p,t}^{\max}$ maximum aggregate demand from consumer group $c$ at time $t$ for water product of state $p$

$ICap_{s,i,p}^{\min}$ minimum capacity to be maintained at source $s$ in inventory $i$ for product of state $p$

$ICap_{s,i,p}^{\max}$ maximum capacity to be maintained at source $s$ in inventory $i$ for product of state $p$

$ICap_{s,i,p}^{buffer}$ buffer capacity limit at source s in inventory $i$ for product of state $p$

$ICap_{s,i,p,t}^{\text{Target}}$ target capacity at time $t$ for product of state $p$ in inventory $i$ at source $s$

$NT$ end time period of the planning horizon

$Op_{s,t}^{ini}$ captures operational state of treatment facility s at the start of the planning horizon

$Q_{s,i,p}^{ini}$ quantity of water product of state $p$ available in inventory $i$ of source $s$ at initial time period of planning horizon

$RS_{r,s}$ assumes value equal to 1 if source $s$ is suitable to transport water in region $r$

$RVP_{r,v,p}$ assumes value equal to 1 if tanker $v$ in region $r$ is compatible to supply product of state $p$

$SC_{s,c}$ assumes value equal to 1 if source $s$ is suitable to supply water to consumer $c$

$SIP_{s,i,p}$ assumes value equal to 1 for suitability of product of state $p$ with inventory $i$ at source $s$

$SMax_s$ Maximum groundwater extraction limit (KL/hour) from source $s$

$SS_{s,s'}$ Suitability of supplying raw water from source $s$ to treatment plant $s'$

$SSP_{s,s',p}$ assumes value 1 if source $s$ is suitable to supply product of state $p$ to treatment plant $s'$

$SSPV_{s,s',p,v}$ assumes value equal to 1 if vehicle type *v* is compatible to supply product *p* from source *s* to treatment plant *s'*

$STpt_s$ throughput of water treatment plant *s* to produce treated water

$STy_s$ indicates source type (FW/GW/TF)

$SP_{s,p}$ assumes value equal to 1 if source *s* is suitable to supply water of state *p*

$T_{s,v}^{Disf}$ water disinfection time for a vehicle type *v* from freshwater source *s*

$T_c^{Distb}$ Distribution time to consumers in a consumer group *c*

$T_s^{DT}$ treatment plant downtime period

$T_{s,v}^{Prep}$ preparation time for vehicle type *v* at source *s*

$T_{s,c,p}^{Transit}$ transit time for transportation of final product water state from source *s* to consumer *c* (i.e. summation of preparation, disinfection and one way travel time)

$T_{s,s',p}^{RWTransit}$ transit time for transportation of raw water state from source *s* to treatment plant *s'* (i.e. summation of preparation and one way travel time)

$T_{s,c,p,v}^{Travel}$ one way travel time required to deliver product state *p* from source *s* to consumer *c* in vehicle type *v*

$T_s^{UT}$ treatment plant uptime period

$TE_t$ captures the start time of a time slot *t*

$TS_t$ captures the end time of a time slot *t*

$TrCost^{Distb}_{s,c,p,v}$ cost of transportation of unit quantity of water in tanker vehicle type *v* for distributing treated water from source *s* to consumer *c*

$TrCost^{RW\,supply}_{s,s',p,v}$ cost of transportation of unit quantity of water in tanker vehicle type *v* for supplying raw water from source *s* to treatment plant *s'*

$TVCost_{s,i,p}$ penalty cost for violating the target capacity of fnal product state *p* in inventory *i* at source *s*

$VA_{r,v,p}$ number of tankers vehicles of type *v* available in region *r* suitable for transporting water product of state *p*

$VExCost_{v,p}$ penalty cost to purchase vehicle type *v* for transporting product of state *p*

$VQ_v$ capacity of tanker vehicle type *v*

$\beta_{s,p}$ fraction of percentage recovery of permeate stream (treated water) from RO process at treatment plant *s*

**Binary Decision Variables (prefix y indicates all the binary variables)**

$yOp_{s,t}$ assumes value equal to 1 if treatment plant is operating at time period *t* to produce treated water

$yPSl_{s,p,t}$ assumes value equal to 1 when product state *p* is selected to be produced from treatment in treatment plant *s'*

**Continuous Decision Variables (prefix x indicates all the continuous variables)**

$xBCV_{s,i,p,t}$ quantity of violation from buffer capacity limit of water product of state *p* in inventory *i* of source *s* at time period *t* $xCDistb_{s,c,p,v}$ quantity of water product of state *p* supplied from source *s* in tanker *v* for distribution to consumer *c*

$xDeCon_{s,c,p,t}$ quantity of water product of state *p* contributed by source *s* at time *t* to supply demand of consumer *c*

$xPDl_{s,c,p,v,t}$ quantity of water product of state *p* delivered to consumer *c* from source *s* in tanker vehicle *v* at time period *t*

$xQ_{s,i,p,t}$ quantity of water product of state *p* available in inventory *i* at time period *t* at source *s*

$xRw_{s,s',p,v,t}$ quantity of raw water supplied from source *s* to treatment plant *s'* in tanker vehicle type *v* at time period *t*

$xSDn_{s,t}$ assumes value equal to 1 if treatment plant *s'* is not operational at time period *t* to treat raw water

$xSSupl_{s,s',p,t}$ quantity of water product of state *p* supplied from source *s* to treatment plant *s'* at time period *t*

$xSUp_{s,t}$ assumes value equal to 1 if treatment plant *s'* is operational at time period *t* to treat raw water

$xTV_{s,i,p,t}$ quantity of violation from target limit of water product *p* in inventory *i* of source *s* at time period *t*

$xTV^{+}_{s,i,p,t}$ quantity of water product of state *p* in inventory *i* of source *s* at time period *t* which is positive violation from target value

$xTV^{-}_{s,i,p,t}$ quantity of water product of state *p* in inventory *i* of source *s* at time period *t* which is negative violation from target value

$xVExQ_{r,v,p}$ extra capacity of tanker vehicle type *v* required in region *r* for transporting product *p*

$xVSSupl_{s,s',p,v}$ quantity of water product of state *p* supplied in tanker vehicle type *v* from source *s* to treatment facility *s'*

**For Table of Contents Only**

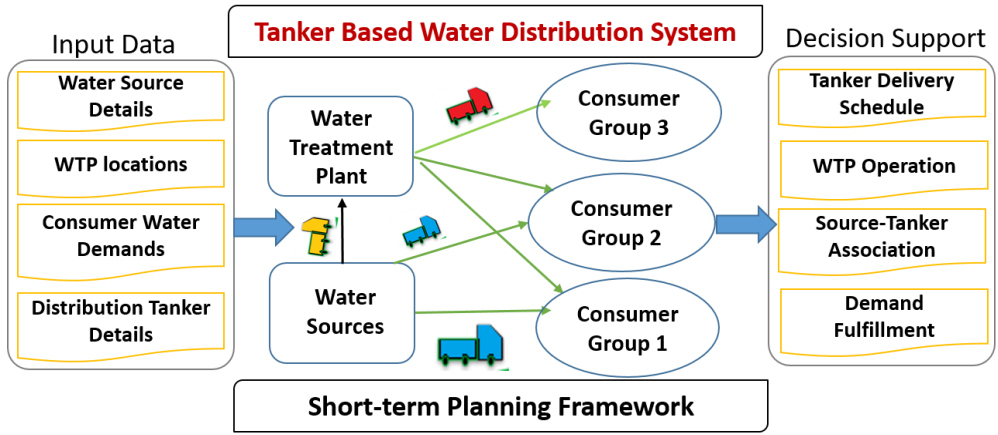